\documentclass[11pt,a4paper]{article}
\usepackage[utf8]{inputenc}
\usepackage[english]{babel}
\usepackage{orcidlink}
\usepackage{graphicx}
\usepackage[left=2cm,right=2cm,top=2cm,bottom=2cm]{geometry}

\usepackage{amsfonts}
\usepackage{amssymb}

\usepackage{bm}
\usepackage{bbm}
\usepackage{mathbbol}
\usepackage{mathrsfs}
\usepackage{amsmath,amsthm,amssymb} 
\usepackage{pifont}

\usepackage{hyperref}
\usepackage{url}

\usepackage{textcomp}

\usepackage{graphicx}
\usepackage{subfigure}

\usepackage{cite}

\usepackage[all,pdf]{xy}

\usepackage{booktabs}
\usepackage{multirow}
\usepackage{multicol}
\usepackage{stfloats}
\usepackage{diagbox}

\newcommand{\blue}[1]{\textcolor{blue}{#1}}
\newcommand{\red}[1]{\textcolor{red}{#1}}

\DeclareMathOperator{\dif}{d}
\DeclareMathOperator{\Dif}{\mathscr{D}}

\newcommand{\mfloor}[1]{ \left\lfloor {#1} \right\rfloor }
\newcommand{\mpair}[2]{ \left\langle {#1}, {#2} \right\rangle}

\newcommand{\me} {\mathrm{e}}

\newcommand{\fracpde}[2]{\frac{\partial {#1}}{\partial {#2}}} 

\newcommand{\fracpdemix}[3]{\frac{\partial^2 {#1}}{\partial {#2} \partial {#3}}}

\newcommand{\set}[1]{\left\{ #1 \right\}}
\newcommand{\abs}[1]{\left| #1 \right|}

\newcommand{\inv}[1]{#1^{-1}}

\newcommand{\ES}[3]{\mathbb{#1}^{{#2}\times {#3}}}     

\newcommand{\scrd}[2]{{#1}_{\mathrm{#2}}}

\newcommand{\BigO}[1]{\mathcal{O}\left(#1\right)}
\newcommand{\PrErrType}[1]{p^\text{#1}_\text{err}}

\usepackage{cleveref}
\usepackage{algorithm}
\usepackage{algorithmicx}
\usepackage{algpseudocode}
\usepackage{float}

\usepackage{forest}

\usepackage{longtable}
\usepackage{supertabular}
\usepackage{caption}

\newcommand{\Algor}{\textbf{Algorithm}~}

\newcommand{\Tab}{\textbf{Table}~}

\algnewcommand\algorithmicswitch{\textbf{switch}}
\algnewcommand\algorithmiccase{\textbf{case}}
\algnewcommand\algorithmicdefault{\textbf{default}}
\algnewcommand\algorithmicassert{\cpvar{assert}}
\algnewcommand\Assert[1]{\State \algorithmicassert(#1)}%
\algdef{SE}[SWITCH]{Switch}{EndSwitch}[1]{\algorithmicswitch\ #1\ \algorithmicdo}{\algorithmicend\ \algorithmicswitch}%
\algdef{SE}[CASE]{Case}{EndCase}[1]{\algorithmiccase\ #1}{\algorithmicend\ \algorithmiccase}%
\algdef{SE}[DEFAULT]{Default}{EndDefault}[1]{\algorithmicdefault\ #1}{\algorithmicend\ \algorithmicdefault}%

\makeatletter
\renewcommand{\ALG@name}{Algorithm}
\newenvironment{breakablealgorithm}
{
	\begin{center}
		\refstepcounter{algorithm}
		\setlength{\baselineskip}{15pt}
		\renewcommand{\caption}[2][\relax]{
			\hrule height.9pt depth0pt \kern3pt
			{\raggedright\textbf{\ALG@name~\thealgorithm} ##2\par}%
			\ifx\relax##1\relax 
			\addcontentsline{loa}{algorithm}{
				\protect\numberline{\thealgorithm}##2}%
			\else 
			\addcontentsline{loa}{algorithm}{
				\protect\numberline{\thealgorithm}##1}%
			\fi
			\kern2pt\hrule\kern2pt
		}
	}{
		\kern3pt\hrule\relax
	\end{center}
}

\usepackage{caption}
\usepackage{listings}
\newcommand{\cpvar}[1]{\texttt{#1}}

\newcommand{\ProcName}[1]{\textsc{#1}}

\author{Hong-Yan Zhang$^{1,2}$\thanks{Corresponding author, email: hongyan@hainnu.edu.cn}\orcidlink{0000-0002-4400-9133}, Zhi-Qiang Feng$^1$, Haoting Liu$^{2,1}$, Rui-Jia Lin$^3$ and Yu Zhou$^1$\\
\begin{tabular}{lll}
$^1$~\small{\textit{School of Information Science and Technology, Hainan Normal University, Haikou 571158, China}}\\
$^2$~\small{\textit{School of Automation, University of Science and Technology
Beijing, Beijing 100083, China}}\\
$^3$~\small{\textit{Information Network and Data Center, Hainan Normal University, Haikou 571158, China}}
\end{tabular}
}

\title{\bf High Order Expansion Method for Kuiper's $V_n$ Statistic in Goodness-of-fit Test}

\date{May 16, 2025}

\begin{document}
\maketitle

\begin{abstract}
Kuiper's $V_n$ statistic, a measure for comparing the difference of ideal distribution and empirical distribution, is of great significance   in the goodness-of-fit test. However, Kuiper's formulae for computing the cumulative distribution function, false positive probability and the upper tail quantile of $V_n$ can not be applied to the case of small sample capacity $n$ since the approximation error is $\BigO{n^{-1}}$. In this work, our contributions lie in three perspectives: firstly the approximation error is reduced  to $\BigO{n^{-(k+1)/2}}$ where $k$ is the expansion order with the \textit{high order expansion}  for the exponent of differential operator; secondly, a novel high order formula with approximation error $\BigO{n^{-3}}$ is obtained by massive calculations; thirdly, the fixed-point algorithms are designed for solving the Kuiper pair of critical values and upper tail quantiles based on the novel formula. The high order expansion method for Kuiper's $V_n$-statistic is applicable for various applications where there are more than five samples of data. The principles, algorithms and code for the high order expansion method are attractive for the goodness-of-fit test.   
\end{abstract}

\tableofcontents

\newpage

\section{Introduction}

In statistics, for two cumulative distribution functions (CDF) $F_1(x)$ and $F_2(x)$, Kuiper's test statistic \cite{Kuiper1960TestsCR,Stephens1970,Pearson1972} is defined by
\begin{equation}
V = D^+ + D^- = \sup_x [F_1(x) - F_2(x)] + \sup_x[F_2(x) - F_1(x)]
\end{equation}
where the discrepancy statistics 
\begin{equation}
D^+ = \sup_{x\in \mathbb{R}} [F_1(x) - F_2(x)], \quad D^- = \sup_{x\in \mathbb{R}} [F_2(x) - F_1(x)] = -\inf_{x\in \mathbb{R}}[F_1(x) - F_2(x)]
\end{equation}
represent the maximum deviation above and below the two CDFs being 
compared, respectively. The trick with Kuiper's test is to use the quantity $D^+ + D^-$  
as the test statistic instead of $\displaystyle \max(D^+, D^-) = \sup_x\abs{F_1(x) - F_2(x)}$ in Kolmogrov-Smirnov test. 
This small change makes Kuiper's test as sensitive in the tails as at the median and 
also makes it invariant under cyclic transformations of the independent random variables. 
The Kuiper's test is useful for the goodness-of-fit test, also named by distribution fitting test  \cite{KStest2008,KnuthTAOCP2,Press1992KuiperTest,
Koch2020Benfords,Dowd2020ECDF,Lanzante2021}. 

Suppose $n$ is the sample capacity, and let 
\begin{equation}
V_n = D^+_n + D^-_n
\end{equation}
be the Kuiper's $V_n$ statistic of critical value for the $V_n$-test. For more details of how to compute the $V_n$ with data collection and the null hypothesis in Kuiper's test, please see \textbf{Appendix} \ref{appendix-Kuiper}.  Kuiper \cite{Kuiper1960TestsCR} pointed out that  the CDF of $V_n$ can be expressed by
\begin{equation} \label{eq-1st-KuiperTest}
\begin{aligned}
F_{V_n}(v) &= \Pr\set{V_n \le v} = \Pr\set{ V_n\cdot \sqrt{n} \le v\cdot \sqrt{n}} \\
 &= \Pr\set{K_n \le c} =\mathscr{A}(c) +\frac{\mathscr{B}(c)}{\sqrt{n}} +  \BigO{\frac{1}{n}}
 \end{aligned}
\end{equation}
 where 
 \begin{equation} \label{eq-c2v}
 \left\{
 \begin{aligned}
  &K_n = \sqrt{n}\cdot V_n \\
  &c  = v\sqrt{n}
 \end{aligned}
 \right.
\end{equation}
and 
\begin{equation}
\mathscr{A}(c) = 1-\sum^\infty_{j=1}2(4j^2c^2-1)\me^{-2j^2c^2}, \quad
\mathscr{B}(c) = \frac{8}{3}c\sum^\infty_{j=1}j^2(4j^2c^2-3)\me^{-2j^2c^2} 
\end{equation}
are functions of $c\in \mathbb{R}^+$. 

Let
\begin{equation} \label{Kuiper-FPP}
\alpha = \Pr\set{K_n>c} = 1 - \mathscr{A}(c) -\frac{\mathscr{B}(c)}{\sqrt{n}} + \BigO{\frac{1}{n}} 
\end{equation}
be the \textit{upper tail probability} (UTP) for the statistics $K_n$ and $V_n$, which is also named the \textit{false positive probability} (FPP) of rejecting the null hypothesis. Equivalently, we have
\begin{equation}
\alpha = \Pr\set{V_n>v} = 1 - \mathscr{A}(v\sqrt{n}) -\frac{\mathscr{B}(v\sqrt{n})}{\sqrt{n}} + \BigO{\frac{1}{n}}
\end{equation} 
by inserting \eqref{eq-c2v} into \eqref{Kuiper-FPP}. The \textit{upper tail quantiles} (UTQ) of $K_n$ and $V_n$, denoted by $c^\alpha_n$ and $v^\alpha_n$ respectively, are the solutions to the following equations
\begin{equation} \label{eq-def-Kuiper-pair}
\left\{
\begin{aligned}
\alpha &= \Pr\set{K_n>c^\alpha_n} =  - \mathscr{A}(c^\alpha_n) -\frac{\mathscr{B}(c^\alpha_n)}{\sqrt{n}} + \BigO{\frac{1}{n}} \\
v^\alpha_n &= c^\alpha_n/\sqrt{n}
\end{aligned}
\right.
\end{equation}
Given the UTP $\alpha$ and the capacity $n$, the solution $\mpair{c^\alpha_n}{ v^\alpha_n}$  to \eqref{eq-def-Kuiper-pair} is called the Kuiper pair. 
We remark that Kuiper provided the tables for the values of $c^\alpha_n$ instead of $v^\alpha_n$ although the $V_n$ is the statistics for the goodness-of-fit test.
By taking the constant term and the first term of $\mathscr{A}(c)$ and $\mathscr{B}(c)$, namely
\begin{equation} \label{eq-Kuiper-AB-approx}
\mathscr{A}(c) \approx 1-(-2+8c^2)\me^{-2c^2}, \quad 
\mathscr{B}(c) \approx \frac{8}{3}(-3c+4c^3)\me^{-2c^2}
\end{equation}
Kuiper obtained the following approximation for the FPP
\begin{equation} \label{eq-alpha-c}
\alpha =  \Pr\set{K_n > c} \approx \left[-2 + 8c^2 - \frac{1}{\sqrt{n}}\cdot \frac{8}{3}(4c^3-3c) \right]\me^{-2c^2}, \quad c > \frac{6}{5}.
\end{equation} 

Although Kuiper provided the tables for the $c^\alpha_n$  
based on \eqref{eq-alpha-c}, there are still some problems that should be considered seriously:
\begin{itemize}
\item[i)] Kuiper's formula for the CDF of $K_n$ expressed by \eqref{eq-1st-KuiperTest} is a $1/2$-order expansion about the sample capacity $n$ and the approximation error $\BigO{n^{-1}}$ implies that the approximation can not be used for small $n$. 
\item[ii)] All of the terms with $\me^{-2j^2c^2}$ of $\mathscr{A}(c)$ and $\mathscr{B}(c)$ for $j\ge 2$ are omitted in \eqref{eq-Kuiper-AB-approx}, which reduces the numerical accuracy  for computing the Kuiper pair $\mpair{c^\alpha_n}{v^\alpha_n}$ if the term with $j=2$ is ignored;
\item[iii)] Kuiper's tables about $\mpair{c^\alpha_n(k)}{v^\alpha_n(k)}$ are just for some values of $n$ and $\alpha$, thus it is impossible to cover versatile scenarios appearing in applications.  
\end{itemize} 
Recently, Zhang et al. \cite{ZhangHY2023KuiperStatistic} improved Kuiper's approximation \eqref{eq-Kuiper-AB-approx} by adding the terms with  $\me^{-8c^2}$ and designed fixed-point algorithms for solving the pair automatically, which deals with the problems ii) and iii) well. Stephens worked on the problem i) along two different directions. 

The first attempt was finished in 1965 by Stephens by proposing two formulae for the Kuiper's $V_n$-statistic \cite{Stephens1965}, namely
\begin{itemize}
\item UTP formula:
\begin{equation} \label{eq-Stephens-utp}
\Pr\set{V_n\ge v} = \sum^{\mfloor{N(1-v)}}_{t=0}\binom{n}{t}(1-v-td)^{n-t-1}W_t
\end{equation}
with  
\begin{equation}
W_t = y^{t-3}\left[y^3 n - y^2 t\frac{3-2/n}{n}-\frac{t(t-1)(t-2)}{n^2} \right], \quad y = v + \frac{t}{n},
\end{equation}
which is valid if $v\ge \frac{1}{2}$ when $n$ is even, and if $v\ge \frac{1}{2} - \frac{1}{2n}$ when $n$ is odd.
\item CDF formula 
\begin{equation} \label{eq-Stephens-approx}
\Pr\set{V_n\le v} =
\left\{
\begin{array}{ll}
n!\left( v-\cfrac{1}{n}\right)^{n-1}, & \cfrac{1}{n} \le v \le \cfrac{2}{n}\\
\cfrac{(n-1)!}{n^{n-2}(t_2-t_1)}\left[t_2^{n-1}(1-t_1) - t_1^{n-1}(1-t_2) \right], & \cfrac{2}{n}\le v \le \cfrac{3}{n}
\end{array}
\right.
\end{equation}
where $t_1$ and $t_2$ are the roots of the quadratic equation of $t$:
$$
t^2 - (nv-1)t + \frac{1}{2}(nv-2)^2 = 0.
$$
\end{itemize}
However, there are three disadvantages for Stephens' formulae. Firstly,
the overestimation/underestimation problem exists. In 2004, Paltani pointed out that Stehphens' formulae will lead to an overestimation of FPP \cite{Paltani2004}. The Python functions \lstinline|kuiper_FPP(D,N)| and \lstinline|kuiper(data, cdf=lambda x: x, args=())| implemented by Anne M. Archibald is based on the Stephens' CDF formula and Paltani's work, in which the overestimation of FPP is emphasized in the comments for the potential usages \cite{Python-kuiper,Python-kuiper-fpp}.
Secondly, the upper bound of approximation error is missing. Unlike Kuiper's work, there is no explicit expression for the approximation error with the variable $n$. 
Finally, the factorials $n!$ and $(n-1)!$ may lead to the overflow problem in computer program when $n$ is large.

The second attempt by Stephens \cite{Stephens1970} appeared in 1970, in which he proposed the modified Kuiper's statistic $T_n$ and the equation for computing the UTQ, namely 
\begin{equation} \label{eq-stephens-Tn}
\left\{
\begin{aligned}
T_n &= V_n\left(\sqrt{n} + 0.155 + \frac{0.24}{\sqrt{n}}\right) = K_n\left(1 + \frac{0.155}{\sqrt{n}} + \frac{0.24}{n} \right)\\
\alpha &= (-2+8c^2)\me^{-2c^2} \Longrightarrow  \scrd{c}{_{MK}}^\alpha=\lim_{n\to\infty}c^\alpha_n
\end{aligned}
\right.
\end{equation}
As the declaration of Stephens in \cite{Stephens1970}, \eqref{eq-stephens-Tn} is applicable for small sample capacity $n$.  We remark that the upper tail quantile $\scrd{c}{_{MK}}^\alpha$ does not depend on $n$ for the modified Kuiper's statistic $T_n$. Unfortunately,
the $T_n$ and $V_n$ will not coincide for large sample capacity since we have $T_n\sim K_n
$ instead of $T_n\sim V_n$ for sufficiently large $n$. Moreover, it is also a pity that there is a lack of description for the correctness and derivation of $T_n$ according to Stephens's work \cite{Stephens1970}.

Unlike the attempts by Stephens, our exploration is along the way of \textit{generalizing the  Kuiper's $\BigO{n^{-1}}$ approximation in \eqref{eq-1st-KuiperTest} to the high order counterpart to get consistent and more accurate equation for solving the Kuiper pair $\mpair{c^\alpha_n}{v^\alpha_n}$ which is applicable for both small and large sample capacity $n$}. There are three key points for our contributions:   
\begin{itemize}
\item[i)]  the approximation error for the Kuiper's formula is reduced  to $\BigO{n^{-(k+1)/2}}$ from $\BigO{n^{-1}}$ where $k$ is the expansion order with the \textit{high order expansion} (HOE) for the exponent of differential operator; 
\item[ii)] a novel high order formula with approximation error $\BigO{n^{-3}}$ is obtained by massive calculations, which can be applied to the scenarios with both large and small sample capacity $n$; 
\item[iii)] the fixed-point algorithms are designed for solving the Kuiper pair of critical values and UTQ based on the novel formula.
\end{itemize}

The rest contents of this paper are organized as follows: the high order expansion for Kuiper's $V_n$-statistic was coped with in Section \ref{sec-hoe-cdf}; 
 the computational method for solving Kuiper pairs with HOE formula is presented in Section \ref{sec-comp-UTQ-hoe}; the fixed-point algorithms for solving the Kuiper pair is dealt with in Section \ref{sec-algorithms-hoe}; the numerical results is given in Section \ref{sec-Num-res-code}; and the conclusion is summarized in Section \ref{sec-conclusion}.

\section{High Order Expansion Formula for Kuiper's $V_n$-Statistic} 
\label{sec-hoe-cdf}

\subsection{High Order Generalization of Kuiper's Probability Formula}

For the purpose of presenting the formula \eqref{eq-1st-KuiperTest}, Kuiper \cite{Darling1957,Kuiper1960TestsCR} introduced the auxiliary function
\begin{equation}
\Phi(a,b) = \sum^\infty_{j=-\infty} \left[ \me^{-2j^2(a+b)^2} - \me^{-2(ja+(j-1)b)^2} \right]
\end{equation}
and the generating function 
\begin{equation} \label{eq-Phin-expansion-1st}
\Phi_n(a,b) = \Phi(a,b) + \frac{1}{6\sqrt{n}}\left(\fracpde{}{a} + \fracpde{}{b} \right)\Phi(a,b) + \BigO{\frac{1}{n}}, \quad a,b \in \mathbb{R}
\end{equation}
for the CDF of $K_n = \sqrt{n}\cdot V_n$.  
Let
\begin{equation}
\Dif_n =  \frac{1}{6\sqrt{n}}\left(\fracpde{}{a} + \fracpde{}{b}\right), \quad \Dif_n^i = \frac{1}{6^i n^{i/2}}\left(\fracpde{}{a} + \fracpde{}{b} \right)^i,
\end{equation}
be the partial derivative operator and its powers respectively, we immediately have
\begin{equation} \label{eq-Phin-1st-approx}
\Phi_n(a,b) = \Phi(a,b) +\Dif_n\Phi(a,b) + \BigO{\frac{1}{n}}.
\end{equation}
Obviously, this is a first expansion for the $\Phi_n(a,b)$ in the sense of the power of $\Dif_n$, or a $1/2$-order expansion for the 
$\Phi_n(a,b)$ in the sense of sample capacity of $n$ since the factor $1/\sqrt{n}$ appears in expression of $\Dif_n$. For the applications with small sample capacity, what we  need is a better expression with smaller approximation error, say $\BigO{n^{-(k+1)/2}}$ for $k=2, 3, \cdots$ where $k$ is a parameter of order for the approximation error. Thus it is necessary for us to expand the $\Phi_n(a,b)$ with more terms which reduces the approximation error as well as makes sure  the expansion is convergent without any extra condition when $k\to +\infty$. It is well known in functional analysis that the norm of the operator $\Dif_n$ is infinite. For the candidate generalization of \eqref{eq-Phin-1st-approx} in the form of series, the domain of convergence should also be infinite since the operator $\Dif_n$ is involved. In consequence, it is a natural strategy for us to replace \eqref{eq-Phin-1st-approx} by the following expression for the purpose of generalization    
\begin{equation} \label{eq-Phin-exp}
\begin{aligned}
\Phi_n(a,b)  &= \exp(\Dif_n)\Phi(a,b)
             = \left[1 + \frac{\Dif_n}{1!} + \frac{\Dif_n^2}{2!} +  \cdots + \frac{\Dif_n^k}{k!}\right]\Phi(a,b) + \BigO{n^{-(k+1)/2}}.
\end{aligned}
\end{equation}
For the $i$-th term in \eqref{eq-Phin-exp} , we can set
\begin{equation} \label{eq-Phi-i-component}
\Phi^{(i)}_n(a,b) = \frac{1}{i!} \Dif_n^i \Phi(a,b) = \frac{1}{n^{i/2}}\cdot \frac{1}{i!6^i}\left[\fracpde{}{a} + \fracpde{}{b} \right]^i\Phi(a,b)
\end{equation}
as the $i$-th component. Hence $\Phi_n(a,b)$ is the superposition of the various components, namely
\begin{equation}
\Phi_n(a,b) = \sum^{k}_{i=0} \Phi^{(i)}_n(a,b) + \BigO{n^{-(k+1)/2}}
\end{equation}
It is easy to find that the computation of $\set{\Phi^{(i)}_n(a,b): i = 0, 1, 2, 3, \cdots}$ can be done iteratively since we have
\begin{equation} \label{eq-iter-deriv}
\Phi^{(r+1)}_n(a,b) = \frac{1}{(r+1)!}\Dif_n^{r+1}\Phi(a,b) 
= \frac{1}{r+1}\Dif_n \Phi^{(r)}_n(a,b)
\end{equation}
for $r=0, 1, 2, \cdots$.

Equation \eqref{eq-Phi-i-component} shows that the continuous variables $a, b$ and the discrete variable $n$ can be separated to distinguish the impacts of the sample capacity $n$ and the label $i$ for the component $\Phi^{(i)}_n(a,b)$. Let
\begin{equation}
Q_i(a,b) = \frac{1}{i!6^i}\left[\fracpde{}{a} + \fracpde{}{b} \right]^i\Phi(a,b),
\end{equation}
then we immediately have
\begin{equation}
\Phi_n(a,b) = \sum^k_{i=0}\frac{Q_i(a,b)}{n^{i/2}} + \BigO{n^{-{(k+1)}/2}}.
\end{equation} 

We now follow Kuiper's derivation \cite{Kuiper1960TestsCR} as follows:
\begin{equation}
\begin{aligned}
\Pr\set{\sqrt{n}\cdot V_n\le c}
&=\int^c_{b=0} \left\{\int^{c-b}_{a=0} \fracpdemix{\Phi_n}{a}{b} \dif a\right\}\dif b
=\int^c_{b=0} \left[\fracpde{\Phi_n}{b} \right]\dif b \\
&= \sum^k_{i=0} \frac{1}{n^i/2} \int^c_{b=0}  \left[\fracpde{Q_i(a,b)}{b}\right]_{a=c-b}\dif b + \BigO{n^{-(k+1)/2}}
\end{aligned}
\end{equation}
Let
\begin{equation}
B_i(c) = \int^c_{b=0} \left[\fracpde{Q_i}{b} \right]_{a=c-b}
\dif b, \quad i = 0, 1, 2, \cdots, k
\end{equation}
then we can obtain the following $k$-th order expansion
\begin{equation} \label{eq-hoe-KuiperTest}
\Pr\set{K_n\le c} 
= \Pr\set{\sqrt{n}\cdot V_n\le c}= \sum^k_{i=0}\frac{B_i(c)}{n^{i/2}} + \BigO{n^{-(k+1)/2}},
\end{equation}
for the CDF $F_{K_n}(c)$ about the sample capacity $n$ with the approximation error $\BigO{n^{-(k+1)/2}}$. Consequently, the relation of $v = c/\sqrt{n}$ implies that the CDF of $V_n$ can be expressed by
\begin{equation} \label{eq-Kuiper-HOE-CDF}
\Pr\set{V_n\le v} = \sum^k_{i=0}\frac{B_i(v\sqrt{n})}{n^{i/2}} + \BigO{n^{-(k+1)/2}}
\end{equation}
 By comparing \eqref{eq-Kuiper-HOE-CDF} with \eqref{eq-1st-KuiperTest} and \eqref{eq-Stephens-approx}, it is easy to find that
\begin{itemize}
\item the HOE formula \eqref{eq-Kuiper-HOE-CDF} is a direct generalization of \eqref{eq-1st-KuiperTest}, where more terms about   the powers of sample capacity $n$ are involved;
\item the approximation error is $\BigO{n^{-(k+1)/2}}$ for \eqref{eq-Kuiper-HOE-CDF}, but the counterpart for \eqref{eq-Stephens-approx} is missing;
\item the computation of \eqref{eq-Stephens-approx} involves the factorial $n!$ and $(n-1)!$, which may lead to the overflow problem caused by large $n$; 
\item the impact of the sample capacity $n$ is separated from the function $B_i(v\sqrt{n})$ in \eqref{eq-Kuiper-HOE-CDF}, which reflects the difference between small capacity and large capacity.
\end{itemize}

\subsection{Expressions for $B_i(c)$}

The principle for computing the functions $B_i(c)$ is not complex but the practical calculation is somewhat tedious and awesome although only the elementary knowledge of calculus is required. For $k=5$, we have
\begin{equation} \label{eq-3rd-KuiperTest}
\Pr\set{K_n\le c} 
= \sum^5_{i=0}\frac{B_i(c)}{n^{i/2}} + \BigO{n^{-3}}
\end{equation}
where 
\begin{equation} \label{eq-Bic-series}
\left\{
\begin{aligned}
B_0(c) & = 1-2\sum_{j=1}^{\infty}(4j^2c^2-1)\me^{-2j^2c^2}\\
B_1(c) & = \dfrac{8}{3}\sum^\infty_{j=1} 
cj^2(4c^2j^2-3)\me^{-2c^2j^2}\\
B_2(c) &= -\frac{1}{18} + \frac{1}{9}\sum^\infty_{j=1}\left\{\me^{-2j^2c^2}\left[
4c^2j^2(-16 c^2 j^4 + 24  j^2 + 1) - 12 j^2- 1\right]\right\} \\
B_3(c) &= \frac{16}{81}\sum^\infty_{j=1}cj^2\me^{-2j^2c^2}\left\{
16 c^4 j^6 -40 c^2 j^4 -4 c^2 j^2 + 15 j^2 + 3\right\} \\
B_4(c) &= \frac{1}{648}+\frac{1}{972}\sum^\infty_{j=1}\left\{\me^{-2j^2c^2}[16c^4j^4(-64 c^2 j^{6} + 240 j^4 + 40 j^2 + 1) \right. \\
&\hspace{2.5cm} \left. - 24c^2j^2(120 j^4 + 40 j^2 + 1) + 120 j^2(2j^2+1) + 3]\right\}\\
B_5(c) &= \frac{2}{3645}\sum^\infty_{j=1}\left\{ 
16\me^{-2c^2j^2}[16 c^5 j^6(32 c^2 j^6 -168 j^4 -40 j^2  - 3)
\right. \\
&\hspace{2.5cm} \left. + 40 c^3 j^4  (84 j^4 + 40 j^2 + 3) - 15cj^2(56 j^4 + 40 j^2 + 3)]
\right\}
\end{aligned}
\right.
\end{equation}
It is prosaic that the $B_i(c)$ is a generalization of Kuiper's $\mathscr{A}(c)$ and $\mathscr{B}(c)$ since  $B_0(c) = \mathscr{A}(c)$ and $B_1(c) = \mathscr{B}(c)$. Particularly, for $k=1$, \eqref{eq-hoe-KuiperTest} is the same with \eqref{eq-1st-KuiperTest}. For more details about the computation of $B_i(c)$, please go to the GitHub site 
\textcolor{blue}{\href{https://github.com/GrAbsRD/KuiperVnHOE}{https://github.com/GrAbsRD/KuiperVnHOE}} for the supplementary material. 

Equation \eqref{eq-3rd-KuiperTest} is a third order formula for computing the CDF, which is much more accurate than the first order version obtained by Kuiper.  For $n=6$, the approximation error is about $6^{-3} \approx 4.63\times 10^{-3}$, thus the probability calculated with \eqref{eq-3rd-KuiperTest} must be  accurate enough for solving practical problems in science, technology, engineering and mathematics. In other words, \eqref{eq-3rd-KuiperTest}  can be applied to the scenario when the sample capacity $n$ is relative small. Although we can get more accurate result for  $k>5$ according to \eqref{eq-hoe-KuiperTest} theoretically, it is not necessary to do so due to the satisfactory accuracy for $k=5$ and boring complexity of calculating the $B_i(c)$ for large $i$.

\section{Computation of Upper Tail Quantile}
\label{sec-comp-UTQ-hoe}

The method for computing the Kuiper pair is similar with the fixed-point 
algorithms stated by Chen et al. \cite{ZhangHY2023KuiperStatistic}, but the order $k$ must be introduced to solve the Kuiper pair $\mpair{c^\alpha_n(k)}{v^\alpha_n(k)}$ via the HOE formula \eqref{eq-Kuiper-HOE-CDF}. 

\subsection{Approximation of $B_i(c)$}

The typical feature of $B_i(c)$ is the property of rapid convergence due to the exponential functions $\me^{-2j^2c^2}$ for $j\in \mathbb{N}$. According to the strategy used in \cite{ZhangHY2023KuiperStatistic}, we take the first and second terms involves $\me^{-2c^2}$ and $\me^{-8c^2}$.  In other words, we have
the following approximation
\begin{equation} \label{eq-Bi-n-k}
B_i(c) \approx \sum^{2}_{j=0}B^j_i(c)\me^{-2j^2c^2} = B^0_i(c) + B^1_i(c)\me^{-2c^2} + B^2_i(c)\me^{-8c^2}  
\end{equation}
where the $B^j_i(c)$ denotes the $j$-th coefficient of $\me^{-2j^2c^2}$ in the expression of $B_i(c)$ for $j \in \set{0, 1, 2}$.  
Substituting \eqref{eq-Bic-series} into \eqref{eq-Bi-n-k}, we can obtain the following approximation equations
\begin{equation}
\left\{
\begin{aligned}
B_0(c) &\approx 1 -  (8c^2-2)\me^{-2c^2} - (32c^2-2)\me^{-8c^2}\\
B_1(c) &\approx  \frac{8}{3}c(4c^2-3)\me^{-2c^2} + \frac{32}{3}c(16c^2-3)\me^{-8c^2}\\
B_2(c) &\approx -\frac{1}{18} -\frac{1}{9}[4c^2(16c^2-25)+13]\me^{-2c^2} - \frac{1}{9}[16c^2(256c^2-97)+49]\me^{-8c^2}\\
B_3(c) &\approx \frac{32}{81}c(8c^4-22c^2+9)\me^{-2c^2} + \frac{64}{81}c(1024c^4 - 656c^2 + 63)\me^{-8c^2}\\
B_4(c) &\approx \frac{1}{648} +\frac{1}{972}[16c^4(-64c^2+281)-3864c^2 + 363]\me^{-2c^2}\\
 &\hspace{1cm} +\frac{1}{972}[256c^4(-4096c^2+4001)-199776c^2+2403]\me^{-8c^2}\\
B_5(c)&\approx\frac{32}{3645}[16c^5(32c^2-211)+5080c^3-1485c]\me^{-2c^2}\\
&\hspace{1cm} +\frac{32}{3645}[1024c^5(2048c^2-2851)+964480c^3-63540c]
\me^{-8c^2}
\end{aligned}
\right.
\end{equation}

\subsection{Fixed-point Equation for the Upper Tail Quantile}

For the FFP $\alpha = \Pr\set{\sqrt{n}\cdot V_n > c}=1-\Pr\set{\sqrt{n}\cdot V_n \le c}$ we have
\begin{equation}
\begin{aligned}
1-\alpha &= \sum^k_{i=0}\frac{B_i(c)}{n^{i/2}}+\BigO{n^{-(k+1)/2}} \\
&\approx \sum^k_{i=0}\sum^2_{j=0}\frac{B^j_i(c)}{n^{i/2}}\me^{-2j^2c^2}+\BigO{n^{-(k+1)/2}} \\
&=\sum^2_{j=0}\left[\sum^k_{i=0}\frac{B^j_i(c)}{n^{i/2}}\right]\me^{-2j^2c^2}+\BigO{n^{-(k+1)/2}} \\
&=\sum^k_{i=0}\frac{B^0_i(c)}{n^i/2} + \me^{-2c^2}\sum^k_{i=0}\frac{B^1_i(c)}{n^{i/2}}
+ \me^{-8c^2}\sum^k_{i=0}\frac{B^2_i(c)}{n^{i/2}}+\BigO{n^{-(k+1)/2}}. 
\end{aligned}
\end{equation}
Hence, we approximately have 
\begin{equation} \label{eq-alpha-A0A1A2}
\alpha \approx \left[1-\sum^k_{i=0}\frac{B^0_i(c)}{n^{i/2}}\right] -\me^{-2c^2}\sum^k_{i=0}\frac{B^1_i(c)}{n^{i/2}} - \me^{-8c^2}\sum^k_{i=0}\frac{B^2_i(c)}{n^{i/2}}
\end{equation}
with the error bound $\BigO{n^{-(k+1)/2}}$.
Let
\begin{equation}
\begin{aligned}
A_j(c,n,k) 
&= -\sum^k_{i=0}\frac{B^j_i(c)}{n^{i/2}}\\
&= -\left[B^j_0(c) +\frac{B^j_1(c)}{n^{1/2}} +\frac{B^j_2(c)}{n} + \frac{B^1_3(c)}{n^{3/2}} + \cdots + \frac{B^1_k(c)}{n^{k/2}}\right], 
\end{aligned}
\end{equation}
then the following iterative formulae hold for computing $A^{(k)}_i$, i.e.,
\begin{equation} \label{eq-Aj-iter}
A_j(c, n, r+1) = A_j(c,n, r) - \frac{B_{r+1}^j(c)}{n^{(r+1)/2}}, \quad j= 0, 1, 2
\end{equation} 
for $r = 0, 1, 2, \cdots, k-1$. It should be noted that each $B^0_i(c)$ is an expression about $n$ and does not depend on $c$, thus we can omit the variable $c$ and denote it as $B^0_i$. Therefore, for $j\in \set{0,1,2}$ we have
\begin{equation} \label{eq-def-Ainck}
\left\{
\begin{aligned}
A_0(n,k)   & = -\left[B^0_0 +\frac{B^0_1}{n^{1/2}} +\frac{B^0_2}{n} + \frac{B^0_3}{n^{3/2}} + \cdots + \frac{B^0_k}{n^{k/2}}\right] \\
A_1(c,n,k) &= -\left[B^1_0(c) +\frac{B^1_1(c)}{n^{1/2}} +\frac{B^1_2(c)}{n} + \frac{B^1_3(c)}{n^{3/2}} + \cdots + \frac{B^1_k(c)}{n^{k/2}}\right] \\ 
A_2(c,n,k) &=
 -\left[B^2_0(c) +\frac{B^2_1(c)}{n^{1/2}} +\frac{B^2_2(c)}{n} + \frac{B^2_3(c)}{n^{3/2}} + \cdots +\frac{B^2_k(c)}{n^{k/2}}\right]
\end{aligned}
\right.
\end{equation}
By substituting \eqref{eq-def-Ainck} into \eqref{eq-alpha-A0A1A2}, we can obtain
\begin{equation} \label{eq-alpha-A0A1A2-2}
\alpha = \left[1+A_0(n,k)\right] +A_1(c,n,k)\me^{-2c^2} + A_2(c,n,k) \me^{-8c^2}.
\end{equation}
Multiplying the factor $\me^{2c^2}$ on the two sides of \eqref{eq-alpha-A0A1A2-2}, we immediately have
\begin{equation}
\me^{2c^2}\left\{\alpha - \left[1+A_0(n,k)\right] \right\} = A_1(c,n,k) + A_2(c,n,k)\me^{-6c^2}.
\end{equation}
Equivalently, we have the fixed-point equation 
\begin{equation} \label{eq-alpha-approx-k}
\boxed{2c^2 +\ln \left[\alpha -1 - A_0(n,k) \right]= \ln\left[A_1(c,n,k) + A_2(c,n,k)\me^{-6c^2}\right] }
\end{equation}
for the  UTQ $c^\alpha_n(k)$ for $K_n$ such that 
\begin{equation}
\alpha = \Pr\set{K_n>c^\alpha_n(k)}
\end{equation}
Once the $c^\alpha_n(k)$ is solved, the 
UTQ $v^\alpha_n(k)$ for $V_n$ can be determined by \eqref{eq-c2v}, namely
\begin{equation}
v^\alpha_n(k) = \frac{c^\alpha_n(k)}{\sqrt{n}}.
\end{equation} 

Particularly, for $k=5$ we have the most important result for the practical applications about HOE method for the Kuiper test
\begin{equation} \label{eq-def-A-B-k5}
\left\{
\begin{aligned}
A_0(n,5) =& -1+\frac{1}{18n}-\frac{1}{648n^2}\\
A_1(c,n,5) =&  (8c^2-2)  - \frac{8}{3\sqrt{n}}(4c^3-3c) \\
&+ \frac{1}{9n}(64 c^4 - 100c^2 + 13) \\
& - \frac{32}{81n^{3/2}}( 8 c^5 - 22 c^3 + 9c)\\
&+ \frac{1}{972n^2}(1024 c^6 - 4496 c^4 + 3864 c^2 - 363) \\
&- \frac{32}{3645n^{5/2}}(512 c^7 - 3376 c^5 + 5080 c^3 - 1485 c)
\\
A_2(c, n, 5) =& (32c^2-2) - \frac{32}{3\sqrt{n}}(16c^3-3c)\\
& + \frac{1}{9n}(4096 c^4 - 1552c^2 + 49) \\
& - \frac{64}{81n^{3/2}}( 1024 c^5 - 656 c^3 + 63c) \\
&+ \frac{1}{972n^2}(1048576 c^6 - 1024256 c^4 + 199776 c^2 - 2403)\\
& - \frac{32}{3645n^{5/2}}(2097152 c^7 - 2919424 c^5 + 964480 c^3 - 63540 c)
\end{aligned}
\right.
\end{equation} 

\subsection{Fixed-Point Method for Solving Upper Tail Quantile}

We can obtain the critical value $c^\alpha_n(k)$ and the upper tail quantile $v^\alpha_n(k) = c^\alpha_n(k)/\sqrt{n}$ by solving the fixed-point  of \eqref{eq-alpha-approx-k} with iterative methods.

\subsubsection{Direct Iterative Method}
Let
\begin{equation}
\scrd{f}{ctm}(c,\alpha,n,k) = \sqrt{\frac{\ln\left[A_1(c,n,k) + A_2(c,n,k)\me^{-6c^2}\right] - \ln \left[\alpha - 1 -A_0(n,k) \right]}{2}},
\end{equation}
then $c$ is the fixed-point of 
\begin{equation} \label{eq-fixedpoint-c-direct-iter}
c = \mathcal{A}^{\alpha,k}_n(\scrd{f}{ctm},c)= \scrd{f}{ctm}(c,\alpha,n,k)
\end{equation}
We can solve the fixed-point  of \eqref{eq-alpha-approx-k} with the following iterative formula 
\begin{equation} \label{eq-solve-c-dit}
c_{i+1} = \mathcal{A}^{\alpha,k}_n(\scrd{f}{ctm},c_i)= \scrd{f}{ctm}(c_i,\alpha,n, k), \quad i = 0, 1, 2, 3, \cdots
\end{equation}
with an appropriate initial value $c_0$.

\subsubsection{Newton Iterative Method}
Let
\begin{equation} \label{eq-Ganc-def}
\scrd{f}{nlm}(c, \alpha, n, k)  = 2c^2 +\ln \left[\alpha - 1 -A_0(n,k)\right]- \ln\left[A_1(c,n,k) + A_2(c,n,k)\me^{-6c^2}\right],
\end{equation}
then the Newton's iterative formula for the fixed-point can be  written by
\begin{equation}
c_{i+1} = \mathcal{B}^{\alpha,k}_n(\scrd{f}{nlm}, c_i)
\end{equation}
where
\begin{equation} \label{eq-solve-c-nit}
\mathcal{B}^{\alpha,k}_n(\scrd{f}{nlm}, c) = c - \frac{\scrd{f}{nlm}(c, \alpha, n, k)}{\scrd{f}{nlm}'(c, \alpha, n, k)}
\end{equation}
is the Newton's updating function. 
The mapping $\mathcal{B}^{\alpha,k}_n(\scrd{f}{nlm}, c)$ is a contractive mapping for appropriate domain of $c$.

\section{Fixed-Point Algorithms for Solving Kuiper Test $\mpair{c^\alpha_n(k)}{v^\alpha_n(k)}$}
\label{sec-algorithms-hoe}

The algorithms discussed below are used to solve the fixed-point determined by \eqref{eq-alpha-approx-k} with the direct iterative method based on \eqref{eq-solve-c-dit} and the Newton iterative method based on \eqref{eq-solve-c-nit} respectively. The fundamental principle, unified framework and algorithms for solving the fixed-point of a nonlinear equation are included in \textbf{Appendix} \ref{app-fixedpoint-algor}.

\subsection{Auxiliary Procedures}

In order to reduce the structure complexity of our key algorithms for solving the Kuiper pair $\mpair{c^\alpha_n(k)}{v^\alpha_n(k)}$, it is wise to introduce some auxiliary procedures. The procedures \ProcName{FunA0}, \ProcName{FunA1}  and \ProcName{FunA2} in \Algor \ref{alg-calc-A0}, \Algor \ref{alg-calc-A1} and
\Algor \ref{alg-calc-A2}  are designed for computing $A_0(n,k)$,  $A_1(c, n, k)$ and $A_2(c, n, k)$ respectively. Note that the iterative formula \eqref{eq-Aj-iter} is used for these procedures.

\begin{breakablealgorithm}
\caption{Calculating the value of $A_0(n)$ arising in the updating function in Kuiper's $V_n$-test}
\label{alg-calc-A0}
\begin{algorithmic}[1]
\Require sample capacity $n$, positive order $k\in \mathbb{N}$ with the default value $k=1$. 
\Ensure  the value of $A_0(n, k)$.
\Function{FunA0}{$n, k$}
\State $A_0\gets -1$; 
\If{$k>1$}
\State $A_0 \gets A_0 + 1.0/(18n)$; \Comment{In C/C++, $1/(18*n)$ will be $0$.}
\EndIf
\If{$k>3$}
\State $A_0 \gets A_0 - 1.0/(648n^2)$; \Comment{In C/C++, $1/(648*n*n)$ will be $0$.}
\EndIf
\State \Return $A_0$;
\EndFunction
\end{algorithmic}
\end{breakablealgorithm}

\begin{breakablealgorithm}
\caption{Calculating the value of $A_1(c, n, k)$ arising in the updating function in $V_n$-test}
\label{alg-calc-A1}
\begin{algorithmic}[1]
\Require  positive real variable $c$,  sample capacity $n$, positive order $k\in \set{1,2,3,4,5}$ with the default value $k=1$. 
\Ensure  the value of $A_1(c, n, k)$.
\Function{FunA1}{$c, n, k$}
\State $ A_1 \gets  (8c^2-2)  - 8(4c^3-3c)/(3\sqrt{n}) $;
\If{$k>1$}
\State $A_1 \gets A_1 + (64 c^4 - 100c^2 + 13)/(9n)$;
\EndIf
\If{$k>2$}
\State $A_1\gets A_1 - 32( 8 c^5 - 22 c^3 + 9c)/(81n^{3/2})$;
\EndIf
\If{$k>3$}
\State $A_1 \gets A_1 + (1024 c^6 - 4496 c^4 + 3864 c^2 - 363)/(972n^2)$;
\EndIf
\If{$k>4$}
\State $A_1 \gets A_1 - 32(512 c^7 - 3376 c^5 + 5080 c^3 - 1485 c)/(3645n^{5/2})$;
\EndIf
\State \Return $A_1$;
\EndFunction
\end{algorithmic}
\end{breakablealgorithm}

\begin{breakablealgorithm}
\caption{Calculating the value of $A_2(c, n, k)$ arising in the updating function in $V_n$-test}
\label{alg-calc-A2}
\begin{algorithmic}[1]
\Require  positive real variable $c$,  sample capacity $n$, positive order $k\in \set{1,2,3,4,5}$ with the default value $k=1$. 
\Ensure  the value of $A_2(c, n, k)$.
\Function{FunA2}{$c, n, k$}
\State $ A_2 \gets (32c^2-2) - 32(16c^3-3c)/(3\sqrt{n}) $;
\If{$k>1$}
\State $A_2 \gets A_2 + (4096 c^4 - 1552c^2 + 49)/(9n)$;
\EndIf
\If{$k>2$}
\State $A_2\gets A_2 - 64( 1024 c^5 - 656 c^3 + 63c)/(81n^{3/2})$;
\EndIf
\If{$k>3$}
\State $A_2 \gets A_2 + (1048576 c^6 - 1024256 c^4 + 199776 c^2 - 2403)/(972n^2)$;
\EndIf
\If{$k>4$}
\State $A_2 \gets A_2 - 32(2097152 c^7 - 2919424 c^5 + 964480 c^3 - 63540 c)/(3645n^{5/2}) $;
\EndIf
\State \Return $A_2$;
\EndFunction
\end{algorithmic}
\end{breakablealgorithm}

The procedure \ProcName{FunFnlm} in  \Algor \ref{alg-f-nlm}   is designed for computing the non-linear functions $\scrd{f}{nlm}(c, \alpha, n, k)$ for the $V_n$-test.

\begin{breakablealgorithm}
\caption{Calculating the value of $\scrd{f}{nlm}(c, \alpha, n, k)$ for the updating function in $V_n$-test}
\label{alg-f-nlm}
\begin{algorithmic}[1]
\Require the UTP $\alpha$,  sample capacity $n$, positive  $c$, positive order $k\in \mathbb{N}$ with default value $k=1$.
\Ensure  the value of $\scrd{f}{nlm}(c, \alpha, n, k)$.
\Function{FunFnlm}{$c, \alpha, n, k$}
\State $A_0 \gets \ProcName{FunA0}(n,k)$;
\State $ A_1 \gets \ProcName{FunA1}(c,n,k)$;
\State $ A_2 \gets \ProcName{FunA2}(c,n,k)$;
\State $ u \gets 2c^2 + \ln (\alpha -1-A_0) - \ln\left(A_1 + A_2\cdot \me^{-6c^2}\right)$;
\State \Return $u$;
\EndFunction
\end{algorithmic}
\end{breakablealgorithm}

The procedure \ProcName{FunFctm}  in  \Algor \ref{alg-f-ctm} is designed for computing the non-linear functions $\scrd{f}{ctm}(c, \alpha, n, k)$ for the $V_n$-test.

\begin{breakablealgorithm}
\caption{Calculating the $\scrd{f}{ctm}(c,\alpha,n,k)$ for the updating operator $\mathcal{A}^{\alpha,k}_n(f, c)$ in $V_n$-test }
\label{alg-f-ctm}
\begin{algorithmic}[1]
\Require  critical value $c$, the UTP $\alpha$, integer $n$,  integer $k$
\Ensure  the value of the direct updating function $\mathcal{A}^{\alpha,k}_n(f, c)$.
\Function{FunFctm1}{$c, \alpha, n, k$}
\State $A_0 \gets \ProcName{FunA0}(n,k)$;
\State $ A_1 \gets \ProcName{FunA1}(c,n,k)$;
\State $ A_2 \gets \ProcName{FunA2}(c,n,k)$;
\State $ y \gets\sqrt{\left(\ln(A_1 + A_2 \cdot \me^{-6c^2}) -\ln (\alpha - 1 -A_0) \right)/2}$; 
\State \Return $y$;
\EndFunction
\end{algorithmic}
\end{breakablealgorithm}

\subsection{Updating Mapping } 

The updating mapping is essential for the iterative scheme for solving the fixed-point. The high order procedure \ProcName{UpdateMethodDirect}
in \Algor \ref{alg-Direct-updator} is designed for computing the contractive mapping $\mathcal{A}^{\alpha,k}_n(f, c)$  for the $V_n$-test . Similarly, 
The high order procedure \ProcName{UpdateMethodNewton}
in \Algor \ref{alg-Newton-updator} is designed for computing the contractive mapping $\mathcal{B}^{\alpha,k}_n(f, c)$  for the $V_n$-test with the Newton iterative method.

\begin{breakablealgorithm}
\caption{Calculating the Contractive Mapping $\mathcal{A}^{\alpha,k}_n(f, c)$}
\label{alg-Direct-updator}
\begin{algorithmic}[1]
\Require  Contractive function $\scrd{f}{ctm}(c, \alpha, n, k)$, critical value $c$, UTP $\alpha$, integer $n$, order $k\in \mathbb{N}$ with default value $k=1$.
\Ensure  the value of the direct updating function $\mathcal{A}^{\alpha,k}_n(\scrd{f}{ctm}, c)$. 
\Function{UpdateMethodDirect}{$\scrd{f}{ctm}, c, \alpha, n, k$}
\State $ u \gets \scrd{f}{ctm}(c, \alpha, n, k)$; 
\State \Return $u$;
\EndFunction
\end{algorithmic}
\end{breakablealgorithm}

\begin{breakablealgorithm}
\caption{Calculating the Newton's updating functions $\mathcal{B}^{\alpha,k}_n(\scrd{f}{nlm}, c)$}
\label{alg-Newton-updator}
\begin{algorithmic}[1]
\Require Nonlinear function $\scrd{f}{nlm}(c, \alpha, n, k)$, critical value $c$, UTP $\alpha$, integer $n$, order $k\in \mathbb{N}$ with the default  value $k=1$. 
\Ensure  The value of the Newton's updating function $\mathcal{B}^{\alpha,k}_n(\scrd{f}{nlm},c)$. 
\Function{UpdateMethodNewton}{$\scrd{f}{nlm}, c, \alpha, n,  k$}
\State $ h \gets 10^{-5}$; 
\State $\cpvar{slope} \gets  \left(\scrd{f}{nlm}(c+h, \alpha, n, k) - \scrd{f}{nlm}(c, \alpha, n, k)\right)/h$; \Comment{Calculate  $\scrd{f}{nlm}'(c, \alpha, n, k)$}
\State $\scrd{c}{new} \gets c - \scrd{f}{nlm}(c, \alpha, n, k)/\cpvar{slope}$;
\State \Return $\scrd{c}{new}$;
\EndFunction
\end{algorithmic}
\end{breakablealgorithm}

\subsection{Algorithms for Solving the Kuiper Pair for $V_n$-test}

The procedure \ProcName{KuiperPairSolver} 
listed in \Algor \ref{alg-KuiperPairSolver} provides a unified framework for solving the Kuiper pair $\mpair{c^\alpha_n(k)}{v^\alpha_n(k)}$ with the direct or Newton's iterative method.

\begin{breakablealgorithm}
\caption{Fixed-point iterative algorithm for solving the Kuiper pair $\mpair{c^\alpha_n(k)}{v^\alpha_n(k)}$ in $V_n$-test with the
 direct/Newton iterative method}
\label{alg-KuiperPairSolver}
\begin{algorithmic}[1]
\Require  The sample capacity $n$, UTP $\alpha\in(0, 1)$, guess of the critical value $\scrd{c}{guess}$ with default value $\scrd{c}{guess}=1.8$, integer $\cpvar{method}\in \set{1,2}$ for the Direct/Newton iterative method with default value $\cpvar{method} = 2$ for the Newton's iterative method, order $k\in \mathbb{N}$ 
with the default value $k=1$.  
\Ensure Kuiper pair $\mpair{c^\alpha_n(k)}{v^\alpha_n(k)}$ such that $\alpha = \Pr\set{V_n \ge v^\alpha_n(k)}$ for $V_n$.
\Function{KuiperPairSolver}{$\scrd{c}{guess},\alpha, n, \cpvar{method}, k$}
\If{($\cpvar{method} == 1$)}
   \State $\ProcName{T} \gets \ProcName{UpdateMethodDirect}$; \Comment{using direct iterative method}
   \State $f \gets \ProcName{FunFctm}$; \Comment{Solve the fixed-point by $c = \mathcal{A}^{\alpha,k}_{n}(\scrd{f}{ctm}, c)$ in $V_n$-test}
\Else
   \State $T \gets \ProcName{UpdateMethodNewton}$; \Comment{using Newton's iterative method}
   \State $f\gets \ProcName{FunFnlm}$; \Comment{Solve the fixed-point by $c = \mathcal{B}^{\alpha,k}_{n}(\scrd{f}{nlm}, c)$ in $V_n$-test}
\EndIf

\State $\epsilon \gets 10^{-5}$;
\State $d\gets \ProcName{Distance}$; \Comment{distance function}
\State $c^\alpha_n(k) \gets \ProcName{FixedPointSolver}(T, f, d, \epsilon, \scrd{c}{guess}, \alpha, n, k)$;
\State $v^\alpha_n(k) \gets c^\alpha_n(k)/\sqrt{n}$;
\State \Return $\mpair{c^\alpha_n(k)}{v^\alpha_n(k)}$; \Comment{Kuiper pair for  the $V_n$-test}
\EndFunction
\end{algorithmic}
\end{breakablealgorithm}

Note that when calling the procedure \ProcName{KuiperPairSolver} we can set $\scrd{c}{guess}=1.8$ or compute the $\scrd{c}{guess}$ by calling the procedure \ProcName{GetInitValue}, namely
\begin{center}
$\scrd{c}{guess}\gets \ProcName{GetInitValue}(\ProcName{FunFnlm}, 0.6, 3.0, 0.05, \alpha, n, k)$;
\end{center}

\subsection{Algorithms for Solving Quantiles and Inverse of CDF}

For most applications, we may have no interest in the critical value $c^\alpha_n(k)$ and our emphasis is put on the UTQ $v^\alpha_n(k)$ and \textit{lower tail quanitle} (LTQ) $v_{1-\alpha}^n(k)$ such that $ v^\alpha_n(k)= v_{1-\alpha}^n(k)$ for the  $V_n$-statistic. Moreover, it is enough to choose the Newton's iterative method in order to find the UTQ or LTQ. In \Algor \ref{alg-Kuiper-UTQ},
the procedure \ProcName{KuiperUTQ} is designed to solve the upper tail quantile. We remark that in the implementation of \ProcName{KuiperUTQ}, the procedure \ProcName{FixedPointSolver} is used instead of \ProcName{KuiperPairSolver}. 

\begin{breakablealgorithm}
\caption{Computing the upper tail quantile in Kuiper's $V_n$-test}
\label{alg-Kuiper-UTQ}
\begin{algorithmic}[1]
\Require the UTP $\alpha\in(0, 1)$, capacity $n$ of the samples, order $k\in \mathbb{N}$.  
\Ensure  the UTQ $v^\alpha_n(k)$ for the $V_n$ statistic such that $\alpha = \Pr\set{V_n > v^\alpha_n(k)}$. 
\Function{KuiperUTQ}{$\alpha, n, k$}
\If{$\alpha \ge  0.9999$}
\State \Return $0.0$;
\EndIf
\State $T \gets \ProcName{UpdateMethodNewton}$; \Comment{using Newton's iterative method}
\State $f\gets \ProcName{FunFnlm}$; \Comment{solve the nonlinear eq. $\scrd{f}{nlm}(c, \alpha, n, k) = 0$ in $V_n$-test}
\State $\epsilon \gets 10^{-5}$; \Comment{precision for the fixed-point}
\State $d\gets \ProcName{Distance}$; \Comment{distance function}
\State $\scrd{c}{guess} \gets 1.8$; \Comment{ or $\scrd{c}{guess}\gets \ProcName{GetInitValue}(\ProcName{FunFnlm}, 0.6, 3.0, 0.05, \alpha, n, k)$; }
\State $c^\alpha_n(k) \gets \ProcName{FixedPointSolver}(T, f, d, \epsilon, \scrd{c}{guess}, \alpha, n, k)$;
\State $v^\alpha_n(k) \gets  c^\alpha_n(k) /\sqrt{n}$; 
\State \Return $v^\alpha_n(k)$; 
\EndFunction
\end{algorithmic}
\end{breakablealgorithm}

The LTQ may become more attractive for some applications.  The procedure 
\ProcName{KuiperLTQ} listed in \Algor \ref{alg-Kuiper-LTQ} is based on the procedure
\ProcName{KuiperPairSolver}. An alternative implementation 
of \ProcName{KuiperLTQ} can be done with the equivalence of $v_\alpha^n(k) = v^{1-\alpha}_n(k)$. In other words, we have $\ProcName{KuiperLTQ}(\alpha,n,k) \equiv \ProcName{KuiperUTQ}(1-\alpha, n,k)$ for any $\alpha \in [0, 1]$.

\begin{breakablealgorithm}
\caption{Compute the Lower Tail Quantile $v_\alpha^n(k)$ in Kuiper's $V_n$-test}
\label{alg-Kuiper-LTQ}
\begin{algorithmic}[1]
\Require the LTP $\alpha\in(0, 1)$, capacity $n$ of the samples, order $k$.  
\Ensure  the LTQ $v_\alpha^n(k)$ in  Kuiper's $V_n$-test. 
\Function{KuiperLTQ}{$\alpha, n, k$}
\If{$\alpha \le  0.0001$}
\State \Return $0.0$;
\EndIf
\State $\scrd{c}{guess} \gets 1.8$; \Comment{or $\scrd{c}{guess}\gets \ProcName{GetInitValue}(\ProcName{FunFnlm}, 0.6, 3.0, 0.05, \alpha, n, k)$;}
\State $\cpvar{method}\gets 2$; \Comment{ Newton's iterative method}
\State $\mpair{c}{v}\gets \ProcName{KuiperPairSolver}(\scrd{c}{guess}, 1-\alpha, n, \cpvar{method}, k) $;
\State \Return $v$; \Comment{lower tail quantile, $v = v_{\alpha}^n(k) = v^{1-\alpha}_n(k)$ }
\EndFunction
\end{algorithmic}
\end{breakablealgorithm}

For Kuiper's goodness-of-fit test, it is the 
inverse of the CDF of $V_n$ that must be solved. It is easy to find that
\begin{equation}
v_x^n = \inv{F}_{V_n}(x), \quad \forall x\in [0,1]
\end{equation}
by the definitions of \textit{lower tail probability} (LTP) and CDF. Therefore, the inverse of CDF can be obtained by calling the procedure \ProcName{KuiperLTQ} or \ProcName{KuiperUTQ} directly with the required arguments. 
The procedure \ProcName{KuiperInvCDF} in \Algor \ref{alg-Kuiper-InvCDF} is designed for computing the  $\inv{F}_{V_n}(x)$ for any probability $x\in[0,1]$.

\begin{breakablealgorithm}
\caption{Compute the inverse CDF $\inv{F}_{V_n}(x)$ for Kuiper test with the HOE formula}
\label{alg-Kuiper-InvCDF}
\begin{algorithmic}[1]
\Require Probability $x \in[0, 1]$, sample capacity $n$, order $k\in \mathbb{N}$
\Ensure  The value of $\inv{F}_{V_n}(x)$ of the  Kuiper's $V_n$-test.
\Function{KuiperInvCDF}{$x, n, k$}
\State $ y \gets \ProcName{KuiperUTQ}(1.0 - x, n, k)$; \Comment{$ y \gets \ProcName{KuiperLTQ}(x, n, k)$ is OK}
\State \Return $y$;
\EndFunction
\end{algorithmic}
\end{breakablealgorithm}

\subsection{Initial Value for Iterative Computation}

Usually it is a combination of art and technique for selecting the initial value when solving nonlinear equation with iterative method. For solving the Kuiper pair $\mpair{c^\alpha_n(k)}{v^\alpha_n(k)}$ and the upper/lower tail of probability involved, it is lucky that the initial value $\scrd{c}{gauss} = 1.8$ is a good choice for $n\ge 6$, $1\le k\le 5$ and $\forall \alpha\in (0,1)$ in the sense of science and statistics. 

An alternative method for finding the initial value is the bisection method for solving the nonlinear equation \eqref{eq-Ganc-def}. The \Algor \ref{alg-find-guess} can be used  by setting $f \gets \scrd{f}{nlm}$, $(x, \mu_1, \cdots, \mu_r) \gets (c, \alpha, n, k)$, $a \gets 0.6$, $b \gets 3$ and $h \gets 0.05$.

\section{Numerical Results} 

\label{sec-Num-res-code}

\subsection{Tables for Kuiper Pair $\mpair{c^\alpha_n(k)}{v^\alpha_n(k)}$}

Kuiper provided numerical tables for the sample capacity $n\in \set{10, 20, 30, 40, 100, +\infty}$ and the UTP $\alpha \in \set{0.01, 0.05, 0.10}$ \cite{Kuiper1960TestsCR}. Chen et al \cite{ZhangHY2023KuiperStatistic} presented consistent results and provides more results for the order $k=1$. It is significant for us to compare the variation of Kuiper pair with the expansion order $k$ for the given $n$ and $\alpha$. \Tab \ref{tab-KuiperPair-alpha-001}, \Tab \ref{tab-KuiperPair-alpha-005}, \Tab \ref{tab-KuiperPair-alpha-010},  \Tab \ref{tab-KuiperPair-alpha-015}, \Tab \ref{tab-KuiperPair-alpha-020}, \Tab \ref{tab-KuiperPair-alpha-030} and \Tab \ref{tab-KuiperPair-alpha-040} demonstrate the Kuiper pairs for $\alpha \in \set{0.01, 0.05, 0.10, 0.20, 0.30, 0.40}$ respectively with expansion order $k\in \set{1, 2, 3, 4, 5}$. It is not difficult for us to obtain the following observations based on these tables:
\begin{itemize}
\item For $k=1$, the Kuiper pairs $\mpair{c^\alpha_n(1)}{v^\alpha_n(1)}$ coincide with those provided by Kuiper \cite{Kuiper1960TestsCR}. We remark that this case was discussed by Zhang et al \cite{ZhangHY2023KuiperStatistic}, in which the error about $(\alpha, n, k) = (0.01, 30, 1)$ in Kuiper's work \cite{Kuiper1960TestsCR} was discovered and fixed. 
\item For the given $\alpha$ and $n$, $v^\alpha_n(k)$ increases monotonically for $k = 1, 2, 3$, which implies that the probability of acceptance in Kuiper's goodness-of-fit test will increases if $k$ increases. 
\item The value $v^\alpha_n(k)$ tends to saturate for $3\le k \le 5$ and there are some fluctuation around $k=4$, which implies that the probability of acceptance for $k = 4$ increases may be slightly smaller than its counterpart for $k=3$. 
\item For $n\ge 100$, $\mpair{c^\alpha_n(k)}{v^\alpha_n(k)}$ are approximately equal for $2\le k\le 5$ when $\alpha$ and $n$ are given, which implies that the second order expansion is enough for improving the accuracy of Kuiper pair. 
\item For $n\ge 2000$, the improvement of accuracy achieved by the HOE method is limited, which means that the original version of Kuiper's test is enough and we can consider the sample is of large scale. 
\item For small $\alpha$ and small $n$, the Kuiper pair $\mpair{c^\alpha_n(k)}{v^\alpha_n(k)}$ varies clearly and significantly with the expansion order $k$. 
\end{itemize}

\begin{table*}[htb]
\centering
\caption{Kuiper pair $\mpair{c^\alpha_n(k)}{v^\alpha_n(k)}$ for expansion order $1\le k \le 5$ and UTP $\alpha = 0.01$}
\label{tab-KuiperPair-alpha-001}
\resizebox{\textwidth}{!}{
\begin{tabular}{|c|ccccc|}
\hline
\diagbox{$n$}{\rotatebox{-5}{$\mpair{c^\alpha_n(k)}{v^\alpha_n(k)}$}}{$k$} &  $1$ &  $2$ &  $3$ &  $4$ & $5$ \\ 
\hline
$6     $ & $(1.7573, 0.7174)$& $(2.2430, 0.9157)$& $(2.1942, 0.8958)$& $(2.1957, 0.8964)$& $(2.1918, 0.8948)$\\ 
$7     $ & $(1.7871, 0.6755)$& $(2.1118, 0.7982)$& $(2.0835, 0.7875)$& $(2.0849, 0.7880)$& $(2.0865, 0.7886)$\\ 
$8     $ & $(1.8092, 0.6397)$& $(2.0625, 0.7292)$& $(2.0427, 0.7222)$& $(2.0436, 0.7225)$& $(2.0454, 0.7231)$\\ 
$9     $ & $(1.8264, 0.6088)$& $(2.0355, 0.6785)$& $(2.0205, 0.6735)$& $(2.0211, 0.6737)$& $(2.0227, 0.6742)$\\ 
\blue{$10$} & $(\blue{1.8401}, 0.5819)$& $(2.0186, 0.6383)$& $(2.0067, 0.6346)$& $(2.0071, 0.6347)$& $(2.0084, 0.6351)$\\ 
\blue{$20$} & $(\blue{1.9026}, 0.4254)$& $(1.9752, 0.4417)$& $(1.9721, 0.4410)$& $(1.9722, 0.4410)$& $(1.9724, 0.4410)$\\ 
\blue{$30$} & $(\underline{\red{1.9252}}, 0.3515)$& $(1.9706, 0.3598)$& $(1.9690, 0.3595)$& $(1.9690, 0.3595)$& $(1.9691, 0.3595)$\\ 
\blue{$40$} & $(\blue{1.9374}, 0.3063)$& $(1.9704, 0.3115)$& $(1.9694, 0.3114)$& $(1.9694, 0.3114)$& $(1.9695, 0.3114)$\\ 
\hline
$45    $ & $(1.9418, 0.2895)$& $(1.9707, 0.2938)$& $(1.9699, 0.2937)$& $(1.9699, 0.2937)$& $(1.9700, 0.2937)$\\
$50    $ & $(1.9454, 0.2751)$& $(1.9712, 0.2788)$& $(1.9705, 0.2787)$& $(1.9705, 0.2787)$& $(1.9705, 0.2787)$\\ 
$55    $ & $(1.9484, 0.2627)$& $(1.9717, 0.2659)$& $(1.9711, 0.2658)$& $(1.9711, 0.2658)$& $(1.9711, 0.2658)$\\ 
$60    $ & $(1.9510, 0.2519)$& $(1.9722, 0.2546)$& $(1.9717, 0.2545)$& $(1.9717, 0.2545)$& $(1.9717, 0.2545)$\\ 
$180   $ & $(1.9739, 0.1471)$& $(1.9806, 0.1476)$& $(1.9805, 0.1476)$& $(1.9805, 0.1476)$& $(1.9805, 0.1476)$\\ 
$200   $ & $(1.9754, 0.1397)$& $(1.9814, 0.1401)$& $(1.9814, 0.1401)$& $(1.9814, 0.1401)$& $(1.9814, 0.1401)$\\ 
$250   $ & $(1.9783, 0.1251)$& $(1.9831, 0.1254)$& $(1.9830, 0.1254)$& $(1.9830, 0.1254)$& $(1.9830, 0.1254)$\\ 
$300   $ & $(1.9804, 0.1143)$& $(1.9844, 0.1146)$& $(1.9843, 0.1146)$& $(1.9843, 0.1146)$& $(1.9843, 0.1146)$\\ 
$400   $ & $(1.9833, 0.0992)$& $(1.9863, 0.0993)$& $(1.9862, 0.0993)$& $(1.9862, 0.0993)$& $(1.9862, 0.0993)$\\ 
$500   $ & $(1.9853, 0.0888)$& $(1.9876, 0.0889)$& $(1.9876, 0.0889)$& $(1.9876, 0.0889)$& $(1.9876, 0.0889)$\\ 
$550   $ & $(1.9860, 0.0847)$& $(1.9882, 0.0848)$& $(1.9882, 0.0848)$& $(1.9882, 0.0848)$& $(1.9882, 0.0848)$\\ 
$560   $ & $(1.9862, 0.0839)$& $(1.9883, 0.0840)$& $(1.9883, 0.0840)$& $(1.9883, 0.0840)$& $(1.9883, 0.0840)$\\
$600   $ & $(1.9867, 0.0811)$& $(1.9887, 0.0812)$& $(1.9886, 0.0812)$& $(1.9886, 0.0812)$& $(1.9886, 0.0812)$\\ 
$700   $ & $(1.9878, 0.0751)$& $(1.9895, 0.0752)$& $(1.9895, 0.0752)$& $(1.9895, 0.0752)$& $(1.9895, 0.0752)$\\ 
$800   $ & $(1.9887, 0.0703)$& $(1.9901, 0.0704)$& $(1.9901, 0.0704)$& $(1.9901, 0.0704)$& $(1.9901, 0.0704)$\\ 
$900   $ & $(1.9894, 0.0663)$& $(1.9907, 0.0664)$& $(1.9907, 0.0664)$& $(1.9907, 0.0664)$& $(1.9907, 0.0664)$\\ 
$1000  $ & $(1.9900, 0.0629)$& $(1.9912, 0.0630)$& $(1.9912, 0.0630)$& $(1.9912, 0.0630)$& $(1.9912, 0.0630)$\\ 
$2000  $ & $(1.9933, 0.0446)$& $(1.9939, 0.0446)$& $(1.9939, 0.0446)$& $(1.9939, 0.0446)$& $(1.9939, 0.0446)$\\ 
\blue{$10^6$} & $(\blue{2.0006}, 0.0020)$& $(2.0006, 0.0020)$& $(2.0006, 0.0020)$& $(2.0006, 0.0020)$& $(2.0006, 0.0020)$\\ 
\hline
\end{tabular}
}
\end{table*}

\begin{table*}[htb]
\centering
\caption{Kuiper pair $\mpair{c^\alpha_n(k)}{v^\alpha_n(k)}$ for expansion order $1\le k \le 5$ and UTP $\alpha = 0.05$}
\label{tab-KuiperPair-alpha-005}
\resizebox{\textwidth}{!}{
\begin{tabular}{|c|ccccc|}
\hline
\diagbox{$n$}{\rotatebox{-5}{$\mpair{c^\alpha_n(k)}{v^\alpha_n(k)}$}}{$k$} &  $1$ &  $2$ &  $3$ &  $4$ & $5$ \\ 
\hline
$6     $ & $(1.5490, 0.6324)$& $(1.6529, 0.6748)$& $(1.6461, 0.6720)$& $(1.6446, 0.6714)$& $(1.6516, 0.6742)$\\ 
$7     $ & $(1.5688, 0.5930)$& $(1.6559, 0.6259)$& $(1.6506, 0.6239)$& $(1.6496, 0.6235)$& $(1.6541, 0.6252)$\\ 
$8     $ & $(1.5842, 0.5601)$& $(1.6590, 0.5865)$& $(1.6547, 0.5850)$& $(1.6540, 0.5848)$& $(1.6571, 0.5859)$\\ 
$9     $ & $(1.5965, 0.5322)$& $(1.6619, 0.5540)$& $(1.6584, 0.5528)$& $(1.6578, 0.5526)$& $(1.6601, 0.5534)$\\ 
\blue{$10$} & $(\blue{1.6066}, 0.5080)$& $(1.6647, 0.5264)$& $(1.6617, 0.5255)$& $(1.6613, 0.5253)$& $(1.6630, 0.5259)$\\ 
\blue{$20$} & $(\blue{1.6563}, 0.3704)$& $(1.6833, 0.3764)$& $(1.6823, 0.3762)$& $(1.6822, 0.3762)$& $(1.6825, 0.3762)$\\ 
\blue{$30$} & $(\blue{1.6758}, 0.3060)$& $(1.6932, 0.3091)$& $(1.6927, 0.3090)$& $(1.6927, 0.3090)$& $(1.6928, 0.3091)$\\ 
\blue{$40$} & $(\blue{1.6868}, 0.2667)$& $(1.6996, 0.2687)$& $(1.6992, 0.2687)$& $(1.6992, 0.2687)$& $(1.6993, 0.2687)$\\ 
$50$ & $(1.6940, 0.2396)$& $(1.7041, 0.2410)$& $(1.7039, 0.2410)$& $(1.7038, 0.2410)$& $(1.7039, 0.2410)$\\ 
\blue{$100$} & $(\blue{1.7110}, 0.1711)$& $(1.7159, 0.1716)$& $(1.7158, 0.1716)$& $(1.7158, 0.1716)$& $(1.7158, 0.1716)$\\ 
\blue{$10^6$} & $(\blue{1.7469}, 0.0017)$& $(1.7469, 0.0017)$& $(1.7469, 0.0017)$& $(1.7469, 0.0017)$& $(1.7469, 0.0017)$\\ 
\hline
\end{tabular}
}
\end{table*}

\begin{table*}[htb]
\centering
\caption{Kuiper pair $\mpair{c^\alpha_n(k)}{v^\alpha_n(k)}$ for expansion order $1\le k \le 5$ and UTP $\alpha = 0.10$}
\label{tab-KuiperPair-alpha-010}
\resizebox{\textwidth}{!}{
\begin{tabular}{|c|ccccc|}
\hline
\diagbox{$n$}{\rotatebox{-5}{$\mpair{c^\alpha_n(k)}{v^\alpha_n(k)}$}}{$k$} &  $1$ &  $2$ &  $3$ &  $4$ & $5$ \\ 
\hline
$6     $ & $(1.4390, 0.5875)$& $(1.5026, 0.6134)$& $(1.5025, 0.6134)$& $(1.5006, 0.6126)$& $(1.5052, 0.6145)$\\ 
$7     $ & $(1.4554, 0.5501)$& $(1.5098, 0.5707)$& $(1.5095, 0.5706)$& $(1.5082, 0.5700)$& $(1.5113, 0.5712)$\\ 
$8     $ & $(1.4684, 0.5192)$& $(1.5158, 0.5359)$& $(1.5155, 0.5358)$& $(1.5145, 0.5354)$& $(1.5167, 0.5362)$\\ 
$9     $ & $(1.4789, 0.4930)$& $(1.5209, 0.5070)$& $(1.5205, 0.5068)$& $(1.5198, 0.5066)$& $(1.5214, 0.5071)$\\ 
\blue{$10$} & $(\blue{1.4877}, 0.4704)$& $(1.5253, 0.4824)$& $(1.5249, 0.4822)$& $(1.5243, 0.4820)$& $(1.5255, 0.4824)$\\ 
\blue{$20$} & $(\blue{1.5322}, 0.3426)$& $(1.5505, 0.3467)$& $(1.5502, 0.3466)$& $(1.5501, 0.3466)$& $(1.5503, 0.3467)$\\ 
\blue{$30$} & $(\blue{1.5503}, 0.2830)$& $(1.5623, 0.2852)$& $(1.5621, 0.2852)$& $(1.5621, 0.2852)$& $(1.5621, 0.2852)$\\ 
\blue{$40$} & $(\blue{1.5606}, 0.2468)$& $(1.5695, 0.2482)$& $(1.5694, 0.2481)$& $(1.5694, 0.2481)$& $(1.5694, 0.2481)$\\ 
$50    $ & $(1.5675, 0.2217)$& $(1.5745, 0.2227)$& $(1.5745, 0.2227)$& $(1.5744, 0.2227)$& $(1.5745, 0.2227)$\\ 
\blue{$100$} & $(\blue{1.5838}, 0.1584)$& $(1.5873, 0.1587)$& $(1.5873, 0.1587)$& $(1.5873, 0.1587)$& $(1.5873, 0.1587)$\\ 
\blue{$10^6$} & $(\blue{1.6193}, 0.0016)$& $(1.6193, 0.0016)$& $(1.6193, 0.0016)$& $(1.6193, 0.0016)$& $(1.6193, 0.0016)$\\ 
\hline
\end{tabular}
}
\end{table*}

\begin{table*}[htb]
\centering
\caption{Kuiper pair $\mpair{c^\alpha_n(k)}{v^\alpha_n(k)}$ for expansion order $1\le k \le 5$ and UTP $\alpha = 0.15$}
\label{tab-KuiperPair-alpha-015}
\resizebox{\textwidth}{!}{
\begin{tabular}{|c|ccccc|}
\hline
\diagbox{$n$}{\rotatebox{-5}{$\mpair{c^\alpha_n(k)}{v^\alpha_n(k)}$}}{$k$} &  $1$ &  $2$ &  $3$ &  $4$ & $5$ \\ 
\hline
$6     $ & $(1.3667, 0.5580)$& $(1.4113, 0.5761)$& $(1.4145, 0.5775)$& $(1.4125, 0.5767)$& $(1.4152, 0.5778)$\\ 
$7     $ & $(1.3813, 0.5221)$& $(1.4201, 0.5367)$& $(1.4223, 0.5376)$& $(1.4209, 0.5371)$& $(1.4228, 0.5378)$\\ 
$8     $ & $(1.3930, 0.4925)$& $(1.4272, 0.5046)$& $(1.4288, 0.5052)$& $(1.4278, 0.5048)$& $(1.4292, 0.5053)$\\ 
$9     $ & $(1.4025, 0.4675)$& $(1.4331, 0.4777)$& $(1.4343, 0.4781)$& $(1.4335, 0.4778)$& $(1.4346, 0.4782)$\\ 
$10    $ & $(1.4105, 0.4460)$& $(1.4381, 0.4548)$& $(1.4391, 0.4551)$& $(1.4385, 0.4549)$& $(1.4393, 0.4551)$\\ 
$20    $ & $(1.4518, 0.3246)$& $(1.4658, 0.3278)$& $(1.4660, 0.3278)$& $(1.4658, 0.3278)$& $(1.4660, 0.3278)$\\ 
$30    $ & $(1.4691, 0.2682)$& $(1.4783, 0.2699)$& $(1.4784, 0.2699)$& $(1.4783, 0.2699)$& $(1.4784, 0.2699)$\\ 
$40    $ & $(1.4790, 0.2338)$& $(1.4859, 0.2349)$& $(1.4860, 0.2349)$& $(1.4859, 0.2349)$& $(1.4859, 0.2349)$\\ 
$50    $ & $(1.4856, 0.2101)$& $(1.4912, 0.2109)$& $(1.4912, 0.2109)$& $(1.4911, 0.2109)$& $(1.4912, 0.2109)$\\ 
$100   $ & $(1.5015, 0.1502)$& $(1.5043, 0.1504)$& $(1.5043, 0.1504)$& $(1.5043, 0.1504)$& $(1.5043, 0.1504)$\\ 
$10^6$ & $(1.5366, 0.0015)$& $(1.5366, 0.0015)$& $(1.5366, 0.0015)$& $(1.5366, 0.0015)$& $(1.5366, 0.0015)$\\ 
\hline
\end{tabular}
}
\end{table*}

\begin{table*}[htb]
\centering
\caption{Kuiper pair $\mpair{c^\alpha_n(k)}{v^\alpha_n(k)}$ for expansion order $1\le k \le 5$ and UTP $\alpha = 0.20$}
\label{tab-KuiperPair-alpha-020}
\resizebox{\textwidth}{!}{
\begin{tabular}{|c|ccccc|}
\hline
\diagbox{$n$}{\rotatebox{-5}{$\mpair{c^\alpha_n(k)}{v^\alpha_n(k)}$}}{$k$} &  $1$ &  $2$ &  $3$ &  $4$ & $5$ \\ 
\hline
$6     $ & $(1.3108, 0.5351)$& $(1.3428, 0.5482)$& $(1.3481, 0.5503)$& $(1.3462, 0.5496)$& $(1.3473, 0.5500)$\\ 
$7     $ & $(1.3242, 0.5005)$& $(1.3525, 0.5112)$& $(1.3563, 0.5127)$& $(1.3550, 0.5122)$& $(1.3558, 0.5125)$\\ 
$8     $ & $(1.3349, 0.4719)$& $(1.3603, 0.4809)$& $(1.3632, 0.4820)$& $(1.3622, 0.4816)$& $(1.3628, 0.4818)$\\ 
$9     $ & $(1.3437, 0.4479)$& $(1.3666, 0.4555)$& $(1.3689, 0.4563)$& $(1.3682, 0.4561)$& $(1.3687, 0.4562)$\\ 
$10    $ & $(1.3511, 0.4272)$& $(1.3720, 0.4339)$& $(1.3739, 0.4345)$& $(1.3733, 0.4343)$& $(1.3737, 0.4344)$\\ 
$20    $ & $(1.3901, 0.3108)$& $(1.4011, 0.3133)$& $(1.4016, 0.3134)$& $(1.4015, 0.3134)$& $(1.4016, 0.3134)$\\ 
$30    $ & $(1.4066, 0.2568)$& $(1.4141, 0.2582)$& $(1.4143, 0.2582)$& $(1.4143, 0.2582)$& $(1.4143, 0.2582)$\\ 
$40    $ & $(1.4163, 0.2239)$& $(1.4219, 0.2248)$& $(1.4220, 0.2248)$& $(1.4220, 0.2248)$& $(1.4220, 0.2248)$\\ 
$50    $ & $(1.4227, 0.2012)$& $(1.4272, 0.2018)$& $(1.4273, 0.2018)$& $(1.4273, 0.2018)$& $(1.4273, 0.2018)$\\ 
$100   $ & $(1.4383, 0.1438)$& $(1.4405, 0.1441)$& $(1.4406, 0.1441)$& $(1.4406, 0.1441)$& $(1.4406, 0.1441)$\\ 
$10^6$ & $(1.4730, 0.0015)$& $(1.4730, 0.0015)$& $(1.4730, 0.0015)$& $(1.4730, 0.0015)$& $(1.4730, 0.0015)$\\ 
\hline
\end{tabular}
}
\end{table*}

\begin{table*}[htb]
\centering
\caption{Kuiper pair $\mpair{c^\alpha_n(k)}{v^\alpha_n(k)}$ for expansion order $1\le k \le 5$ and UTP $\alpha = 0.30$}
\label{tab-KuiperPair-alpha-030}
\resizebox{\textwidth}{!}{
\begin{tabular}{|c|ccccc|}
\hline
\diagbox{$n$}{\rotatebox{-5}{$\mpair{c^\alpha_n(k)}{v^\alpha_n(k)}$}}{$k$} &  $1$ &  $2$ &  $3$ &  $4$ & $5$ \\ 
\hline
$6     $ & $(1.2234, 0.4995)$& $(1.2381, 0.5054)$& $(1.2458, 0.5086)$& $(1.2444, 0.5080)$& $(1.2424, 0.5072)$\\ 
$7     $ & $(1.2349, 0.4668)$& $(1.2488, 0.4720)$& $(1.2546, 0.4742)$& $(1.2536, 0.4738)$& $(1.2524, 0.4734)$\\ 
$8     $ & $(1.2443, 0.4399)$& $(1.2572, 0.4445)$& $(1.2618, 0.4461)$& $(1.2610, 0.4458)$& $(1.2603, 0.4456)$\\ 
$9     $ & $(1.2520, 0.4173)$& $(1.2642, 0.4214)$& $(1.2679, 0.4226)$& $(1.2672, 0.4224)$& $(1.2668, 0.4223)$\\ 
$10    $ & $(1.2586, 0.3980)$& $(1.2700, 0.4016)$& $(1.2730, 0.4026)$& $(1.2725, 0.4024)$& $(1.2722, 0.4023)$\\ 
$20    $ & $(1.2940, 0.2893)$& $(1.3008, 0.2909)$& $(1.3017, 0.2911)$& $(1.3015, 0.2910)$& $(1.3015, 0.2910)$\\ 
$30    $ & $(1.3095, 0.2391)$& $(1.3142, 0.2399)$& $(1.3147, 0.2400)$& $(1.3146, 0.2400)$& $(1.3146, 0.2400)$\\ 
$40    $ & $(1.3185, 0.2085)$& $(1.3222, 0.2091)$& $(1.3225, 0.2091)$& $(1.3225, 0.2091)$& $(1.3225, 0.2091)$\\ 
$50    $ & $(1.3247, 0.1873)$& $(1.3277, 0.1878)$& $(1.3279, 0.1878)$& $(1.3278, 0.1878)$& $(1.3278, 0.1878)$\\ 
$100   $ & $(1.3397, 0.1340)$& $(1.3412, 0.1341)$& $(1.3413, 0.1341)$& $(1.3413, 0.1341)$& $(1.3413, 0.1341)$\\ 
$10^6$ & $(1.3740, 0.0014)$& $(1.3740, 0.0014)$& $(1.3740, 0.0014)$& $(1.3740, 0.0014)$& $(1.3740, 0.0014)$\\ 
\hline
\end{tabular}
}
\end{table*}

\begin{table*}[htb]
\centering
\caption{Kuiper pair $\mpair{c^\alpha_n(k)}{v^\alpha_n(k)}$ for expansion order $1\le k \le 5$ and UTP $\alpha = 0.40$}
\label{tab-KuiperPair-alpha-040}
\resizebox{\textwidth}{!}{
\begin{tabular}{|c|ccccc|}
\hline
\diagbox{$n$}{\rotatebox{-5}{$\mpair{c^\alpha_n(k)}{v^\alpha_n(k)}$}}{$k$} &  $1$ &  $2$ &  $3$ &  $4$ & $5$ \\ 
\hline
$6     $ & $(1.1529, 0.4707)$& $(1.1549, 0.4715)$& $(1.1638, 0.4751)$& $(1.1630, 0.4748)$& $(1.1581, 0.4728)$\\ 
$7     $ & $(1.1630, 0.4396)$& $(1.1661, 0.4408)$& $(1.1729, 0.4433)$& $(1.1723, 0.4431)$& $(1.1693, 0.4420)$\\ 
$8     $ & $(1.1712, 0.4141)$& $(1.1750, 0.4154)$& $(1.1804, 0.4173)$& $(1.1799, 0.4171)$& $(1.1779, 0.4164)$\\ 
$9     $ & $(1.1781, 0.3927)$& $(1.1822, 0.3941)$& $(1.1866, 0.3955)$& $(1.1862, 0.3954)$& $(1.1848, 0.3949)$\\ 
$10    $ & $(1.1841, 0.3744)$& $(1.1882, 0.3757)$& $(1.1919, 0.3769)$& $(1.1916, 0.3768)$& $(1.1906, 0.3765)$\\ 
$20    $ & $(1.2166, 0.2720)$& $(1.2200, 0.2728)$& $(1.2211, 0.2731)$& $(1.2210, 0.2730)$& $(1.2209, 0.2730)$\\ 
$30    $ & $(1.2311, 0.2248)$& $(1.2337, 0.2252)$& $(1.2343, 0.2254)$& $(1.2343, 0.2253)$& $(1.2342, 0.2253)$\\ 
$40    $ & $(1.2397, 0.1960)$& $(1.2418, 0.1964)$& $(1.2422, 0.1964)$& $(1.2422, 0.1964)$& $(1.2422, 0.1964)$\\ 
$50    $ & $(1.2456, 0.1762)$& $(1.2474, 0.1764)$& $(1.2476, 0.1764)$& $(1.2476, 0.1764)$& $(1.2476, 0.1764)$\\ 
$100   $ & $(1.2601, 0.1260)$& $(1.2611, 0.1261)$& $(1.2612, 0.1261)$& $(1.2612, 0.1261)$& $(1.2612, 0.1261)$\\ 
$10^6$ & $(1.2939, 0.0013)$& $(1.2939, 0.0013)$& $(1.2939, 0.0013)$& $(1.2939, 0.0013)$& $(1.2939, 0.0013)$\\ 
\hline
\end{tabular}
}
\end{table*}

\subsection{Probability of Type I Errors}


The order of HOE for Kuiper's $V_n$-statistic has significant impacts on the \textit{probability of Type I errors}, which is denoted by $\PrErrType{I}$. In order to demonstrate this impact, we take the standard normal population and sample data to validate the performance of interest. 

Let $X_1,X_2,\cdots,X_{\scrd{n}{rep}} \overset{\textrm{i.i.d.}}{\sim}\mathcal{N}(0,1) $.  We can perform $ \scrd{n}{rep}$ independent Monte Carlo experiments to calculate the value of $\PrErrType{I}$. For the $i$-th Kuiper's goodness-of-fit test with significance level $ \alpha $, we can set the following hypothesis:
\begin{itemize}
\item null hypothesis $\mathcal{H}_0$, namely the data set $\set{x_1, x_2, \cdots, x_{\scrd{n}{rep}}}$ generated form the $X_1,X_2,\cdots,X_{\scrd{n}{rep}}$ satisfy the  distribution $\mathcal{N}(0,1)$;
\item  alternative hypothesis $\mathcal{H}_1$, namely the data set $\set{x_1, x_2, \cdots, x_{\scrd{n}{rep}}}$ does not satisfy the distribution $\mathcal{N}(0,1)$.
\end{itemize}
Let $\mathcal{E}^{\mathrm{I}}_i(k)$ be the indicator variable for the Type I error with $k$-th order method such that
\begin{equation}
\mathcal{E}^{\mathrm{I}}_i(k) = \begin{cases}
1, & \mathcal{H}_0 \text{~is rejected --- Type I error occurs}\\
0, & \mathcal{H}_0 \text{~is accepted --- no  error}\\
\end{cases}
\end{equation}
for $1\le i\le \scrd{n}{rep}$, then the $\PrErrType{I}$ for the $k$-th order method can be calculated by 
\begin{equation*}\label{eq-prob-type-I-error}
\PrErrType{I}(k)=\dfrac{1}{\scrd{n}{rep}}\sum_{i=1}^{\scrd{n}{rep}}\mathcal{E}^{\mathrm{I}}_i(k)
\end{equation*}
according to the large number's law if $\scrd{n}{rep}$ is large enough. In our experiment, we set $\scrd{n}{rep} = 1000$. The simulation result is shown in \Tab \ref{tab-type-I-errors}.   
It is trivial to obtain some observations for the same sample capacity $n$:
\begin{itemize}
\item For $1\le k\le 5$ and $n\ge 6$, we always have $\PrErrType{I}(k) \le \alpha = 0.05$. However, for the Kolmogrov-Smirnov test we have
$\PrErrType{I} > \alpha$ and for $n\in \set{7, 9, 10, 30, 50, 100}$ we can find that $\PrErrType{I}$ is near $\alpha$ although $\PrErrType{I}<\alpha $. Moreover, for $n\in \set{6,7}$ the $\PrErrType{I}$ for the Stephens's test is larger than $\alpha$. 
\item For Kuiper's test, we have $\PrErrType{I}(2)=\PrErrType{I}(3)=\PrErrType{I}(4)=\PrErrType{I}(5)$ due to the Kuiper pairs $\mpair{c^{0.05}_n(k)}{v^{0.05}_n(k)}$ are approximately the same for $2\le k \le 5$ according to \Tab \ref{tab-KuiperPair-alpha-005}.  In other words, \Tab \ref{tab-KuiperPair-alpha-005} and \Tab \ref{tab-type-I-errors} imply that 
the precision of the the Kuiper pair $\mpair{c^{\alpha}_n(k)}{v^{\alpha}_n(k)}$ beyond the third decimal place is irrelevant to the calculation of the probability of Type I errors.
\item For $k\ge 2$, the $\PrErrType{I}(k)$ is much less than the $\PrErrType{I}(1)$, which means HOE method is better than the first order method. 
\item The $\PrErrType{I}(k)$ for the $k$-th order Kuiper's test is much less than that the counterparts of Kolmogrov-Smirnov test and Stephens' test, which implies that Kuiper's test is better in the sense of small sample capacity if the $\PrErrType{I}$ is considered.
\end{itemize} 
Consequently,  for small sample capacity $n$ we should take the HOE method for Kuiper's test with priority when choosing the test method among the Kuiper's test, Kolmogrov-Smirnov test and Stephens' test if the probability of Type I errors is considered.
 
\begin{table*}[htb]
    \centering
    \caption{Type I error probabilities for the expansion order $1\le k \le 5$, significance level $ \alpha=0.05 $, $\hat{q}_t = t/n$,  and repetition times $\scrd{n}{rep} = 1000$ in simulation. The last row is computed by Paltani’s approach with Stephens’ scheme.}
    \label{tab-type-I-errors}
    \resizebox{\textwidth}{!}{
    \begin{tabular}{|c|ccccccccccc|}
        \hline
        \diagbox{$k$}{$ \PrErrType{I} $}{$n$} &  $ 6 $ & $ 7 $ & $ 8 $ & $ 9 $ & $ 10 $ & $ 20 $ & $ 30 $ & $ 40 $ & $ 50 $ & $ 100 $ & $ 180 $ \\ 
        \hline
        $1$ & 
        $0.0040$ & $0.0080$ & $0.0070$ & $0.0100$ & $0.0040$ & $0.0120$ & $0.0120$ & $0.0130$ & $0.0200$ & $0.0320$ & $0.0380$ \\
        $2$ & 
        $0.0010$ & $0.0050$ & $0.0040$ & $0.0060$ & $0.0030$ & $0.0080$ & $0.0110$ & $0.0120$ & $0.0190$ & $0.0290$ & $0.0380$ \\
        $3$ & 
        $0.0010$ & $0.0050$ & ${0.0050}$ & $0.0060$ & $0.0030$ & ${0.0090}$ & $0.0110$ & $0.0120$ & $0.0190$ & $0.0290$ & $0.0380$ \\
        $4$ & 
        $0.0010$ & $0.0050$ & ${0.0050}$ & $0.0060$ & $0.0030$ & ${0.0090}$ & $0.0110$ & $0.0120$ & $0.0190$ & $0.0290$ & $0.0380$ \\
        $5$ & 
        $0.0010$ & $0.0050$ & $0.0040$ & $0.0060$ & $0.0030$ & $0.0080$ & $0.0110$ & $0.0120$ & $0.0190$ & $0.0290$ & $0.0380$ \\
        \hline
        KS-test & 
        $\red{0.0520}$ & $0.0490$ & $\red{0.0520}$ & $0.0480$ & $0.0460$ & $\red{0.0510}$ & $0.0500$ & $\red{0.0590}$ & $0.0490$ & $0.0480$ & $\red{0.0520}$ \\
        \hline
        Stephens & 
        $\red{0.1150}$ & $\red{0.0600}$ & $0.0400$ & $0.0370$ & $0.0300$ & $0.0280$ & $0.0340$ & $0.0270$ & $0.0380$ & $0.0470$ & $0.0480$ \\
        \hline
    \end{tabular}
    }
\end{table*}

\subsection{Discussion}
It should be remarked that the configurations of $\alpha$ and $n$ are correlated in applications. Our suggestions based on the simulation for the parameters $\alpha$ and $n$ are as follows: 
\begin{itemize}
\item for $\alpha = 10^{-2}$, we should set $n\ge 6$; 
\item for $\alpha = 10^{-3}$, we should set $n\ge 56$; 
\item for $\alpha = 10^{-4}$, we should set $n\ge 560$. 
\end{itemize}
Generally, the smaller the $\alpha$ is, the larger the sample capacity $n$ should be. 
The reason for these phenomenons is obvious:  the $\alpha$ means the upper tail significance level or the probability of Type I errors in hypothesis test, hence it is unwise to set tiny $\alpha$ when the sample capacity $n$ is small in the sense of science and statistics. In other words, only large sample capacity $n$ permits small $\alpha$ for the probability of Type I errors.

\section{Conclusion}
\label{sec-conclusion}

Kuiper's formula for computing the CDF, FFP and UTQ is limited for small sample capacity $n$ since its accuracy is controlled by $\BigO{n^{-1}}$. 
With the help of high order expansion of $\Phi_n(a,b) = \exp(\Dif_n)\Phi(a,b)$, the CDF of $V_n$ can be generalized in order to obtain a small upper bound $\BigO{n^{-(k+1)/2}}$ where $k\in \mathbb{N}$ is the expansion order. For $k=5$, the novel formula \eqref{eq-3rd-KuiperTest} and its equivalent forms are obtained through massive calculations to compute the CDF, FFP and UTQ with the approximation error $\BigO{n^{-3}}$, which can be applied directly for the scenarios with small capacity $n$. 

Unlike the Stephens's work \cite{Stephens1965} on the similar objective of improving the accuracy of computing CDF and FFP, our HOE method is more natural and simple. Stephens's modified Kuiper test $T_n$ is not compatible with the large sample capacity due to the asymptotic behavior  $T_n\sim K_n$ for sufficiently large $n$. By comparison, we always have $\displaystyle \lim_{n\to \infty} \Pr\set{K_n\le c } = B_0(c)$,
which leads to the formulae for the large sample capacity naturally. 

With the help of the numerical results for the Kuiper pairs $\mpair{c^\alpha_n(k)}{v^\alpha_n(k)}$, the potential users can set  the order $k$ appropriately according to his/her practical problem. Our recommendation for the choice of the expansion order is $k = 5$ for any $\alpha\in (0, 1)$.

For the convenience of potential users such as researchers, engineers, statisticians and college students, the fixed-point algorithms are also designed and explained briefly. We believe that the pseudo-code for the algorithms and C code for the implementation of the algorithms based on the HOE method are valuable for the readers and users of interest.

\section*{Supplementary information}

The supplementary material is provided for deriving the formulae of $\set{B_i(c): 0\le i\le 5}$, and it  is  also released in the website: \textcolor{blue}{\href{https://github.com/GrAbsRD/KuiperVnHOE}{https://github.com/GrAbsRD/KuiperVnHOE}}.  

\section*{Acknowledgments}

This work was supported in part by the National Natural Science Foundation of China under grant numbers 62167003 and 62373042, in part by the Hainan Provincial Natural Science Foundation of China under grant number 720RC616, in part by the Research Project on Education and Teaching Reform in Higher Education System of Hainan Province under grant numbers Hnjg2025ZD-28 and Hnjg2023ZD-26,  in part by the Specific Research Fund of the Innovation Platform for Academicians of Hainan Province, in part by the Hainan Province Key R \& D Program Project under grant number WSJK2024MS234, in part by the Guangdong Basic and Applied Basic Research Foundation under grant number 2023A1515010275, and in part by the Foundation of National Key Laboratory of Human Factors Engineering under grant number HFNKL2023WW11.

\section*{Declarations}

\subsection*{Conflict of interest}

The authors declare that there is no conflict of interest.

\subsection*{Code availability}

The code for the implementations of the algorithms discussed in this paper can be downloaded from the following GitHub website 
\begin{center}
\textcolor{blue}{\href{https://github.com/GrAbsRD/KuiperVnHOE}{https://github.com/GrAbsRD/KuiperVnHOE}}
\end{center}

\appendix

\section{Goodness-of-fit Test with Kuiper's $V_n$ statistic } \label{appendix-Kuiper}

\subsection{Problem and Assumption}

For the finite $n$ i.i.d samples $X_1, X_2, \cdots, X_n$, a fundamental problem is: whether the samples $\set{X_t: 1\le t\le n}$ satisfy the given continuous distribution of a population $X$ with the distribution function $F_X(x)$? Here we assume that the samples are identically independent variables. 
\subsection{Hypothesis}

For the continuous population $X$ with the CDF $F_X(x)$ and sample set $\set{X_t: 1\le t\le n}$, the null hypothesis $\mathcal{H}_0$ is: the samples $\set{X_t: 1\le t\le n}$ satisfy the distribution $F_X(x)$. 

\subsection{Calculation of $V_n$ with Sample Data}

The Kuiper's $V_n$ statistic is the sum of $D^+_n$ and $D^-_n$, which can be computed with the sample set $\set{X_t: 1\le t\le n}$ directly. For the sample set $\set{X_t: 1\le t\le n}$ of population $X$, the \textit{empirical distribution function} (EDF) will be
\begin{equation}
\hat{F}_X(x)= \frac{1}{n}\sum^{n}_{t=1} \mathbb{1}[X_t \le x]
\end{equation}
where
\begin{equation}
\mathbb{1}[X_t \le x] = \begin{cases}
1, & X_t \le x \\
0, & X_t > x
\end{cases}
\end{equation}
is the indicator function. 
Suppose the $X_{(t)}$ is the  order statistics  such that 
$$X_{(1)} < X_{(2)} < \cdots < X_{(k)} < \cdots  <  X_{(n-1)} < X_{(n)},$$
then the EDF can be calculated easily for the sorted data, that is, 
\begin{equation}
\displaystyle \hat{q}_t = \hat{F}_X\left(X_{(t)}\right) = \hat{q}^{(0)}_t = \frac{1}{n}\sum^{n}_{i=1} \mathbb{1}[X_i \le X_{(t)}] = \frac{t}{n}.
\end{equation}
Let  
\begin{equation}
q_t = F_{X}\left(X_{(t)}\right)
\end{equation} 
then 
\begin{equation}
\left\{
\begin{aligned}
D^+_n &= \max_{1\le t\le n} \left[\hat{F}_X\left(X_{(t)}\right) - F_X\left(X_{(t)}\right)\right] = \max_{1\le t\le n} \left(\hat{q}_t - q_t \right) \\
D^-_n &= \max_{1\le t\le n} \left[F_X\left(X_{(t)}\right) - \hat{F}_X\left(X_{(t)}\right) \right]=  \max_{1\le t\le n} \left(q_t- \hat{q}_t \right)
\end{aligned}
\right.
\end{equation}
Thus we can obtain the value of Kuiper's $V_n$-test statistic
\begin{equation*} 
V_n = D^+_n + D^-_n = \max_{1\le t\le n} \left(\hat{q}_t - q_t \right) +  \max_{1\le t\le n} \left(q_t- \hat{q}_t \right)
\end{equation*}
with the date samples $\set{X_t: 1\le t\le n}$. 

Please note that we can sort the sample data $\set{X_t: 1\le t\le n}$ by a sorting algorithm (such as the quick sorting with computational complexity $\BigO{n\log_2n}$) in an increasing order, thus the $t$-th sorted value is the $X_{(t)}$ of interest. We also remark that the disadvantage for the $\hat{q}_t = t/n$ is that $\hat{q}_n \equiv 1$ when $t=n$.
There  are several choices for computing the empirical probability
$\hat{q}_t$ listed in \Tab \ref{tab-qt} for $1\le t\le n$ and different expressions for $\hat{q}_t$ will lead to different performances for the Kuiper's goodness-of-fit test.
\begin{table}[h]
\centering
\caption{Formula for computing the empirical probability} \label{tab-qt}
\begin{tabular}{cll}
\hline
Scheme & Expression for $\hat{q}_t$ &  Remark\\
\hline
0 &  $\hat{q}^{(0)}_t= \cfrac{t}{n}$ & Used for the computation of $D^+_n$ by Stephens \cite{Stephens1970}\\
1 &  $\hat{q}^{(1)}_t= \cfrac{t-1}{n}$ & Used for the computation of $D^-_n$ by Stephens \cite{Stephens1970}\\
2 &  $\hat{q}^{(2)}_t= \cfrac{t-0.5}{n}$ & \\
3 &  $\hat{q}^{(3)}_t= \cfrac{t}{n+1}$ & \\
4 &  $\hat{q}^{(4)}_t= \cfrac{t-0.375}{n+0.250}$ & ISO standard 5479:1997, Chinese standard GB/T 4882-2001  \\
\hline
\end{tabular}
\end{table}
We remark that the scheme \#4 is recommended by the Chinese standard GB/T 4882-2001 \cite{GB-T4882-2001} or the equivalent international standard ISO 5479:1997. 
Historically, a strategy of mixing was taken \cite{Stephens1970} as follows
\begin{equation} \label{eq-Vn-stephens}
\begin{split}
V_n &= \max_{1\le t\le n} \left(\hat{q}^{(0)}_t - q_t \right) +  \max_{1\le t\le n} \left(q_t- \hat{q}^{(1)}_t \right)\\
&= \max_{1\le t\le n} \left(\cfrac{t}{n} - q_t \right) +  \max_{1\le t\le n} \left(q_t-\cfrac{t-1}{n} \right)
\end{split}
\end{equation}

\section{Methods for Solving Fixed-point}

\label{app-fixedpoint-algor}

\subsection{Principle and Unified theoretic Framework}

For a nonlinear equation
\begin{equation}
f(x) = 0, \quad x\in [a,b],
\end{equation}
usually we can convert it into an equivalent fixed-point equation with the form
\begin{equation}
x = T(f, x), \quad x\in[a,b]
\end{equation}
where $f(\cdot)$ is a mapping and $T(\cdot, \cdot)$ is contractive operator. There are two typical categories of fixed-point equation:
\begin{itemize}
\item[(i)] Direct iterative scheme
      \begin{equation}
         x_{i+1} = T(\scrd{f}{ctm}(x_i), x_i) = \scrd{f}{ctm}(x_i)
      \end{equation}
      which depends on the function $\scrd{f}{ctm}$.
\item[(ii)] Newton's iterative scheme
      \begin{equation}
      x_{i+1} = T(\scrd{f}{nlm}(x_i), x_i) = x_i - \frac{\scrd{f}{nlm}(x_i)}{\scrd{f}{nlm}'(x_i)},
      \end{equation}
      which depends not only on the equivalent nonlinear mapping  
$\scrd{f}{nlm}(x)$ from
\begin{equation}
\phi(x) = 0 \Longleftrightarrow \scrd{f}{nlm}(x) = 0
\end{equation}
but also on its derivative $\scrd{f}{nlm}'(x)$.
\end{itemize}
We remark that the strings ``ctm" and ``nlm" in the subscripts  means \underline{c}on\underline{t}ractive \underline{m}apping and \underline{n}on\underline{l}inear \underline{m}apping respectively.

The stopping condition for the iterative method is satisfied by the Cauchy's criteria for the convergent sequence $\set{x_i: i = 0, 1, 2, 3, \cdots}$
\begin{equation}
d(x_{i+1}, x_i) < \epsilon
\end{equation}
where $\epsilon$ denotes the precision and $d(\cdot,\cdot)$ is a metric for the difference of the current value of $x$, say $\scrd{x}{guess}$ and the its updated version, say $\scrd{x}{improve}$. The simplest choice of $d(\cdot, \cdot)$ is the absolute value function for $x\in [a, b] \in \mathbb{R}$, i.e.,
\begin{equation}
d(x_i, x_{i+1}) = d(x_{i+1}, x_i) = \abs{x_{i+1} - x_i}
\end{equation}

For practical problems arising in various fields, there may be some extra parameters in the function $\phi$, say $\phi(x, \alpha, n)$ where $\alpha$ is a real number, $n$ is an integer and the order of $x, \alpha, n$ is not important\footnote{In other words, both $\phi(x,\alpha,n)$ and $\phi(x,n,\alpha)$ are feasible and correct for computer programming.}. In general, we can use the notations $\phi(x, \mu_1, \cdots, \mu_r), f(x, \mu_1 \cdots \mu_r)$ and $T(f, x, \mu_1, \cdots, \mu_r)$ to represent the scenario when there are some extra parameters. We remark that
the derivative of $f(x,\mu_1, \cdots, m_r)$ with respect to the variable $x$ should be calculated by
\begin{equation}
     f'(x, \mu_1, \cdots \mu_r) = \lim_{\Delta x\to 0} \frac{f(x+\Delta x, \mu_1, \cdots, \mu_r) -f(x,\mu_1, \cdots, \mu_r)}{\Delta x}
\end{equation}
for the Newton's iterative method.

\subsection{Iterative Algorithms}

We now give the pseudo-code for the fixed-point algorithms with the concepts of high order function and function object, please see  \Algor \ref{alg-abstract-fixed-point}. Note that the order for the list of arguments can be set according to the programmer's preferences. 

\begin{breakablealgorithm}
\caption{Unified Framework for Solving the Fixed-Point of $x = \mathcal{T}^*_\star(f, x)= T(f, x, \mu_1, \cdots, \mu_r)$}\label{alg-abstract-fixed-point}
\begin{algorithmic}[1]
\Require Contractive mapping  $T$ as the updator which is a high order function, function object $f$ for the $\scrd{f}{nlm}$ or $\scrd{f}{ctm}$, function object $d$ for the distance $d(x_i, x_{i+1})$, precision $\epsilon$, initial value $\scrd{x}{guess}$ and extra parameters $\mu_1, \cdots, \mu_r$ with the same or different data types.
\Ensure Fixed-point $x$ such that $x = \mathcal{T}^*_\star(f, x) = T(f, x, \mu_1, \cdots, \mu_r)$
\Function{FixedPointSolver}{$T, f, d, \epsilon, \scrd{x}{guess},  \mu_1, \cdots, \mu_r$}
\State $\scrd{x}{improve}\gets T(f, \scrd{x}{guess}, \mu_1, \cdots, \mu_r)$;
\While{$d(\scrd{x}{improve}, \scrd{x}{guess}) \ge \epsilon$}
\State $\scrd{x}{guess} \gets \scrd{x}{improve}$;
\State $\scrd{x}{improve}\gets T(f, \scrd{x}{guess}, \mu_1, \cdots, \mu_r)$;
\EndWhile
\State \Return $\scrd{x}{improve}$;
\EndFunction
\end{algorithmic}
\end{breakablealgorithm}

In the sense of programming language and discrete mathematics, $f$ is an ordinary (first order) function, $T$ is a second order function and $\ProcName{FixedPointSolve}$ is a third order function. 

\begin{breakablealgorithm}
\caption{Calulate the distance of $x$ and $y$}
\label{alg-calc-dist}
\begin{algorithmic}[1]
\Require $x, y\in \mathbb{R}$
\Ensure The distance of $x$ and $y$, i.e., $d(x,y)$ 
\Function{Distance}{$x,y$}
\State $\cpvar{dist} \gets \abs{x -y}$; \Comment{if $x, y\in \ES{R}{n}{1}$, we can use the $\ell^p$-norm}
\State \Return $\cpvar{dist}$;
\EndFunction
\end{algorithmic}
\end{breakablealgorithm}

\subsection{Setting Initial Value With Bisection Method}

The procedure \ProcName{GetInitValue} listed in \Algor \ref{alg-find-guess} is designed for finding an appropriate initial value $\scrd{x}{guess}\in (a, b)$ for root of $f(x, \mu_1, \cdots, \mu_r) =0$ where $f$ is continuous on $(a, b)$ and has a unique root in $(a,b)$. 

\begin{breakablealgorithm}
\caption{Finding the initial value for solving the root of the continuous function $f(x, \mu_1, \cdots, \mu_r) = 0$ which has a unique root in the interval $[a,b]$. } \label{alg-find-guess}
\begin{algorithmic}[1]
\Require Function object $f$, real number $a$ and real number $b$ such that $a<b$, possible extra parameters $\mu_1, \cdots, \mu_r$, precision $h$ with default value $h = 0.05$. 
\Ensure $\scrd{x}{guess}\in(a,b)$ such that it is a good initial value for solving $f(x, \mu_1, \cdots, \mu_r)=0$;
\Function{GetInitValue}{$f, a, b, h, \mu_1, \cdots, \mu_r$}
\State $\Delta \gets \abs{a-b}$; 
\State $\scrd{x}{guess}\gets (a+b)/2$;
\While{$\Delta > h$}
\If{$f(\scrd{x}{guess}, \mu_1, \cdots, \mu_r)\cdot f(a, \mu_1, \cdots, \mu_r) > 0$}
\State $ a\gets \scrd{x}{guess}$;
\Else
\State $b\gets \scrd{x}{guess}$;
\EndIf
\State $\Delta \gets \Delta/2$;
\State $\scrd{x}{guess}\gets (a+b)/2$;
\EndWhile
\State \Return $\scrd{x}{guess}$;
\EndFunction
\end{algorithmic}
\end{breakablealgorithm}


\section*{How to cite this work?}
\begin{itemize}
\item Plain text\\
H.-Y.Zhang, Z.-Q. Feng, H.-T. Liu, R.-J. Liu and Y. Zhou.  High order expansion method for Kuiper’s statistic in goodness-of-fit test. Statistics and Computing, 2025, vol.35, article number 91, pages: 1--19. \url{https://doi.org/10.1007/s11222-025-10623-9}
\item BibTeX
\begin{verbatim}
@Article{ZhangHY2025KuiperHOE,
author = {H.-Y. Zhang and Z.-Q. Feng and H.-T. Liu and R.-J. Lin and Y. Zhou},
title  = {High order expansion method for Kuiper's $V_n$ statistic in 
          goodness-of-fit test},
journal = {Statistics and Computing},
year   = {2025},
volume = {35},
pages  = {e91:1-19},
note   = {\url{https://doi.org/10.1007/s11222-025-10623-9}, Arxiv:2310.19576},
}
\end{verbatim}
\end{itemize}

\end{document}